\documentclass[12pt]{article}

\usepackage{lineno,hyperref}
\usepackage{amsfonts}
\usepackage{amsmath}
\usepackage{amssymb,epsfig}
\usepackage{algorithm}
\usepackage{algorithmic}
\usepackage{ntheorem}
\modulolinenumbers[5]
\newtheorem{thm}{Theorem}

\newtheorem{prop}{Proposition}

\newtheorem*{proof}{Proof}

\def\kld{{\mathcal D}_{\mathrm{KL}}}

\let\<=\langle
\let\>=\rangle

\begin{document}

%\shorttitle{AM-VEnKF}    

%\shortauthors{L. Wen, J. Li}  

\title{Affine-Mapping based Variational Ensemble Kalman Filter}  

%\tnotemark[<tnote number>] 

%\tnotetext[<tnote number>]{<tnote text>} 


\author{Linjie Wen\thanks{The School of Mathematical Sciences, Shanghai Jiao Tong University, Shanghai 200240, China}\, and Jinglai Li\thanks{The School of Mathematics, University of Birmingham, Birmingham B15 2TT, UK}}
\date{}
\maketitle

% Main text
\begin{abstract}
We propose an affine-mapping based variational Ensemble Kalman filter for sequential Bayesian filtering problems with generic observation models. 
Specifically, the proposed method is formulated as to construct an affine mapping from the prior ensemble to the posterior one, 
and the affine mapping is computed via a variational Bayesian formulation, i.e., 
by minimizing the Kullback-Leibler divergence  between the transformed distribution through the affine mapping and the actual posterior. 
Some theoretical properties of resulting optimization problem are studied and 
a gradient descent scheme is proposed to solve the resulting optimization problem.
With numerical examples we demonstrate that the method has competitive performance against existing methods. 
\end{abstract}

\section{Introduction}
The ensemble Kalman filter (EnKF)~\cite{evensen2009data,evensen2003ensemble} is one of the most popular tools for sequential data assimilation, thanks to its computational efficiency and flexibility~\cite{houtekamer1998data,whitaker2002ensemble,evensen2003ensemble}. 
Simply put, at each time step EnKF approximates the prior, the likelihood and the posterior by Gaussian distributions. Such a Gaussian approximation allows an affine update that maps the prior ensemble to the posterior one. This Gaussian approximation and the resulting affine update are the key that enables EnKF to handle large-scale problems with 
a relatively small number of ensembles. 
In the conventional EnKF,  it is required that the observation model is Gaussian-linear, which means that the observation operator is linear
and the noise is additive Gaussian. However, in many real-world applications, neither of these two requirements is satisfied. When the actual observation model is not Gaussian-linear, %it is still treated as such a model, if one wants to use EnKF (see Section~\ref{sec:enkf}), which 
the EnKF method may suffer from substantial estimation error, which is discussed in details in Section~\ref{sec:enkf}.
% Considerably attention has been devoted to extending EKF to the nonlinear and/or non-Gaussian observation model,
% where the majority of such EKF extensions impose certain restrictions on the observation model. 
% For example, the methods in \cite{4982682,arulampalam2002a,reif1999stochastic} can handle the situation where the observation operator is nonlinear while the observation noise is still Gaussian. 
% On the other hand, a number of works study some specific non-Gaussian observation noise, just to name a few~\cite{han2008evaluation,izanloo2016kalman,kitagawa1996monte}.
To the end, it is of practical importance to develop methods that can better deal with {generic observation models}
than EnKF, while retaining the computational advantage (i.e., using a small ensemble size) of it. 

A notable example of such methods is the nonlinear ensemble adjustment filter (NLEAF)~\cite{lei2011moment},
which involves a correction scheme: the posterior moments are calculated with importance sampling and the ensembles are then corrected accordingly. Other methods that can be applied to such problems include \cite{anderson2003local, anderson2001ensemble, houtekamer2001sequential,li2018trimmed,ba2018two} (some of them may need certain modifications), just to name a few. In this work we focus on the EnKF type of methods that can use a small number of ensembles in high dimensional problems, and methods involving full Monte Carlo sampling such as the particle filter (PF)~\cite{arulampalam2002a,doucet2009tutorial} are not in our scope. It is also worth noting that a class of methods combine EnKF and PF to alleviate the estimation bias induced by the non-Gaussianity~(e.g., \cite{stordal2011bridging,frei2013bridging}), and typically the EnKF part in such methods still requires a Gaussian-linear observation model (or to be treated as such a model). 

The main purpose of this work is to provide an alternative framework to implement EnKF for arbitrary observation models. Specifically, the proposed method formulates the EnKF update as to construct an affine mapping from the prior to the posterior and such an affine mapping is computed in variational Bayesian framework~\cite{mackay2003information}. That is, we seek the affine mapping minimizing the Kullback-Leibler divergence (KLD) between the ``transformed'' prior distribution and the posterior. We note here that a similar formulation has been used in the variational (ensemble) Kalman filter~\cite{auvinen2010variational,solonen2012variational}. The difference is however, the variational (ensemble) Kalman filter methods mentioned above still rely on the linear-Gaussian observation model, where the variational formulation, combined with a BFGS scheme, is used to avoid the inversion and storage of very large matrices, while in our work the variational formulation is used to compute the optimal affine mapping for generic observation models. 
%We also note that 

It can be seen that this affine mapping based variational EnKF (VEnKF) reduces to the standard EnKF when the observation model is Gaussian-linear, and as such it is a natural generalization of the standard EnKF to generic observation models. Also, by design the obtained affine mapping is \emph{optimal} under the variational (minimal KLD) principle. We also present a numerical scheme based on gradient descent algorithm to solve the resulting optimization problem, and with numerical examples we demonstrate that the method has competitive performance against several existing methods. Finally we emphasize that, though the proposed method can perform well for generic observation models, it requires the same assumption as the standard EnKF, i.e., the posterior distributions should not deviate significantly from Gaussian. 

The rest of the work is organized as follows. In Section~\ref{sec:formulation} we provide a generic formulation of the sequential Bayesian filtering problem. In Section~\ref{sec:lmekf} we present the proposed affine mapping based variational EnKF. Numerical examples are provided in Section~\ref{sec:examples} to demonstrate the performance of the proposed method and finally some closing remarks are offered in Section~\ref{sec:closing}.

\section{Problem Formulation} \label{sec:formulation}
\subsection{Hidden Markov Model}
We start with the  hidden Markov model (HMM), which is a generic formulation for data assimilation problems~\cite{doucet2009tutorial}. Specifically let $\{x_t\}_{t\geq 0}$ and $\{y_t\}_{t\geq0}$ be two discrete-time stochastic processes, taking values from continuous state spaces $\mathcal{X}$ and $\mathcal{Y}$ respectively. 
Throughout this work we assume that $\mathcal{X}=\mathbb{R}^{n_x}$ and  $\mathcal{Y}=\mathbb{R}^{n_y}$. The HMM model assumes that the pair $\{x_t,y_t\}$ has the following property, 
\begin{subequations}\label{e:hmm}
	\begin{eqnarray}
	x_{t} |x_{1:t-1},y_{1:t-1} &\sim& \pi(x_t|x_{t-1}),\quad x_0\sim \pi(x_0), \label{e:hmm1}\\
	y_t|x_{1:t},y_{1:t-1} &\sim& \pi(y_t|x_t), \label{e:hmm2}
	%\sim f_t(\cdot|x_{t-1}), \label{eq:state}
	\end{eqnarray}
\end{subequations} 
where for simplicity we assume that the probability density functions (PDF) of all the distributions exist and $\pi(\cdot)$ is used as a generic notation of a PDF whose actual meaning is specified by its arguments. 

\begin{figure}[h]
	%\begin{minipage}[htbp]{0.5\linewidth}
	\centerline{\includegraphics[width=0.5\textwidth]{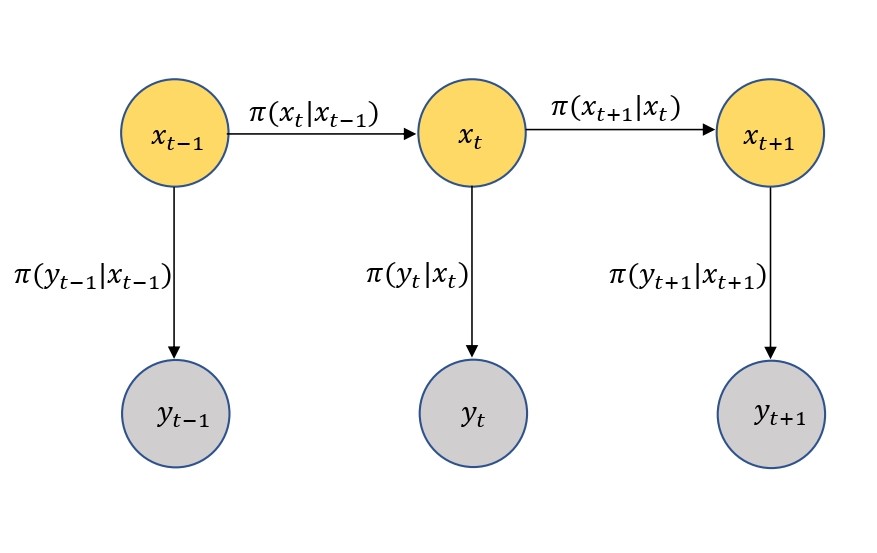}}
	\caption{A schematic illustration of the Hidden Markov Model.}
	\label{fig:hmm}
	
\end{figure}

In the HMM formulation, $\{x_t\}$ and $\{y_t\}$ are  known respectively as the hidden and the observed states, and a schematic illustration of HMM is shown in Fig.~\ref{fig:hmm}. This framework represents many practical problems of interest~\cite{fine1998hierarchical,krogh2001predicting,beal2002infinite}, where one makes observations of $\{y_t\}_{t\geq0}$ and wants to estimate the hidden states $\{x_t\}_{t\geq0}$ therefrom. A typically example of HMM is the following stochastic discrete-time dynamical system:
\begin{subequations}
	\begin{eqnarray}
	x_t &=& F_t(x_{t-1},\alpha_t),\quad x_0\sim \pi(x_0),\label{e:model}\\
	y_t &=& G_t(x_t,\beta_t),\label{e:obs}
	\end{eqnarray}
\end{subequations} 
where $\alpha_t\sim \pi^\alpha_t(\cdot)$ and $\beta_t\sim\pi^\beta_t(\cdot)$ are random variables representing respectively the model error and the observation noise at time $t$. {In many real-world applications
	such as numerical weather prediction~\cite{bauer2015quiet}, Eq.~(\ref{e:model}), which represents the underlying physical model, is computationally intensive, while Eq.~(\ref{e:obs}), describing the observation model, is often available analytically and therefore easy to evaluate. It follows that, in such problems, 1) one can  only afford a small number of particles in the filtering, 2) Eq.~(\ref{e:model}) accounts for the vast majority of the computational cost.} All our numerical examples are described in this form and further details can be found in Section~\ref{sec:examples}.

%	And $f_t(\cdot|x_{t-1})$ denotes the transform probability density from $x_{t-1}$, and the probability measure is denoted by $dx_t$. Secondly, We want to estimate the discrete process $\{x_t\}$ but only have the observation process $\{y_t\}\in\mathcal{Y}$ (the probability measure is denoted by $dy_t$), that is 
%	\begin{equation}
%			y_t|x_t   \sim g_t(\cdot|x_t),\quad t=1,\ldots,T\label{eq:obs}
%	\end{equation}
%	where $T$ is a positive integer, and the stochastic process $\{y_t\}$ is conditional independence on the process $\{x_t\}$. The pair $\{x_t,y_t\}$ composed of Eq.~(\ref{eq:state}) and Eq.~(\ref{eq:obs}) is called (infinite) Hidden Markov model, $f_t(\cdot|x_{t-1})$ and $g_t(\cdot|x_t)$ are known conditional distributions, and $x_t$ and $y_t$ are respectively the hidden and the observed states.
%	
%	
%	We note here that, in many real-world applications, Eq.~(\ref{eq:state}) involves simulating complicated dynamical models, while 
%	Eq.~(\ref{eq:obs}) is much easier to evaluate or simulate, 
%	and therefore Eq.~(\ref{eq:state}) contributes dominantly to the total computational burden. 
%	We assume this is the case in problems that we are interested in. 
%and the prior distribution of initial state $x_0$ is $p(x_0)$. 

\subsection{Recursive Bayesian Filtering}
Recursive Bayesian filtering~\cite{chen2003bayesian} is a popular framework to estimate the hidden states in a HMM, and it aims to compute the condition distribution $\pi(x_t|y_{1:t})$ for $t=1,2,\ldots$ recursively. In what follows we discuss how the recursive Bayesian filtering proceeds. 

First applying the Bayes' formula, we obtain 
\begin{equation}\pi(x_t|y_{1:t}) =\frac{\pi(y_t|x_t,y_{1:t-1})\pi(x_t|y_{1:t-1}) }{\pi(y_{t}|y_{1:t-1})}, \label{e:bayes1}
\end{equation}
where $\pi(y_t|y_{1:t-1})$ is the normalization constant that often does not need to be evaluated in practice. From Eq.~(\ref{e:hmm2}) we know that $y_t$ is independent of $y_{t-1}$ conditionally on $x_t$, and thus Eq.~(\ref{e:bayes1}) becomes 
\begin{equation}\pi(x_t|y_{1:t}) =\frac{\pi(y_t|x_t)\pi(x_t|y_{1:t-1}) }{\pi(y_{t}|y_{1:t-1})}. \label{e:bayes2}
\end{equation}

The condition distribution $\pi(x_t|y_{1:t-1})$ can be expressed as 
\begin{equation}
\pi(x_t|y_{1:t-1})  = \int \pi(x_t|x_{t-1}, y_{1:t-1})\pi(x_{t-1}|y_{1:t-1})dx_{t-1},
\end{equation}
and again thanks to the property of the HMM in Eq.~(\ref{e:hmm}),  we have,
\begin{equation}
\pi(x_t|y_{1:t-1})  = \int \pi(x_t|x_{t-1})\pi(x_{t-1}|y_{1:t-1})dx_{t-1}, \label{e:prior}
\end{equation}
where $\pi(x_{t-1}|y_{1:t-1})$ is the posterior distribution at the previous step $t-1$.

As a result the recursive Bayesian filtering performs the following two steps in each iteration: 
\begin{itemize}
	\item Prediction step: the prior density $\pi(x_t|y_{1:t-1})$ is determined via Eq.~(\ref{e:prior}),
	\item Update step: the posterior density $\pi(x_t|y_{1:t})$  is computed via Eq.~(\ref{e:bayes2}).
\end{itemize}

The recursive Bayesian filtering provides a generic framework for sequentially computing the conditional distribution $\pi(x_t|y_{1:t})$ as the iteration proceeds. In practice, the analytical expressions for the posterior $\pi(x_t|y_{1:t})$ or the prior $\pi(x_t|y_{1:t-1})$ usually can not be obtained, and therefore these distributions have to be represented numerically, for example, by an ensemble of particles. 
%			\begin{figure}[h]
%			%\begin{minipage}[htbp]{0.5\linewidth}
%			\centerline{\includegraphics[width=0.5\textwidth]{figs_LMENKF/RBF}}
%			
%			\caption{An illustration of recursive Bayesian filtering.}
%			\label{fig:ruf}
%		
%		\end{figure}
\section{Affine mapping based VEnKF}\label{sec:lmekf}
We describe the affine-mapping based VEnKF (AM-VEnKF) algorithm in this section.

\subsection{Formulation of the affine-mapping based VEnKF }\label{sec:lmkf}
We first consider the update step: namely suppose that the prior distribution $\pi(x_t|y_{1:t-1})$ is obtained, and we want to compute the posterior $\pi(x_t|y_{1:t})$. 

We start with a brief introduction to the transport map based methods for computing the posterior distribution~\cite{el2012bayesian}, where the main idea is to construct a mapping which pushes the prior distribution into the posterior.  Namely suppose $\tilde{x}_t$ follows the prior distribution  $\pi(\cdot|y_{1:t-1})$, and one aims to construct a bijective mapping {$T: \mathcal{X}\rightarrow \mathcal{X}$}, 
% 	{\color{red}$T: \mathcal{X}\rightarrow \mathcal{X}$,}
such that $x_t=T(\tilde{x}_t)$ follows the posterior distribution $\pi(\cdot|y_{1:t})$. In reality, it is often impossible to exactly push the prior into the posterior $\pi(\cdot|y_{1:t})$, and in this case an approximate approach can be used. That is, let $\pi_T(\cdot)$ be the distribution of $x_t = T(\tilde{x}_t)$ 
where  $\tilde{x}_t\sim \pi(\cdot|y_{1:t-1})$ and we seek a mapping $T\in \mathcal{H}$ where $\mathcal{H}$ is a given function space, so that $\pi_T(\cdot)$ is ``closest'' to the actual posterior $\pi(\cdot|y_{1:t})$ in terms of certain measure of distance between two distributions. 

In practice, the KLD, which (for any two distributions $\pi_1$ and $\pi_2$) is defined as, 
\begin{equation}
\kld(\pi_1,\pi_2) = \int \log\left[ \frac{\pi_1(x)}{\pi_2(x)}\right] \pi_1(x) dx,
\end{equation}
is often used for such a distance measure. 
That is, we find a mapping $T$ by solving the following minimization problem,
\begin{equation}
\min_{T\in\mathcal{H}} \kld(\pi_T,\pi(x_t|y_{1:t})), \label{e:minkld}
\end{equation} 
which can be understood as a variational Bayes formulation. 
%To actually find the mapping, we need to address two key issues. The first is that,  
In practice, the prior distribution $\pi(\tilde{x}_t|y_{1:t-1})$ is usually not analytically available, and in particular they are represented by an ensemble of particles. %As is in the standard EnKF, we approximate the prior with a Gaussian distribution. 
As is in the standard EnKF, we  estimate a Gaussian approximation of the prior distribution $\pi(\tilde{x}_t|y_{1:t-1})$ from the ensemble. Namely, given an ensemble $\{\tilde{x}_t^m\}_{m=1}^M$ drawn from the prior distribution $\hat{\pi}(\tilde{x}_t|y_{1:t-1})$, we construct an approximate prior  $\hat{\pi}(\cdot|y_{1:t-1}) = N(\tilde{\mu}_t,\tilde{\Sigma}_t)$, with 
\begin{equation}\label{eq:emp}
\tilde{\mu}_t=\frac{1}{M}\sum_{m=1}^M\tilde{x}_t^m,\\
\tilde{\Sigma}_t=\frac{1}{M-1}\sum\limits_{m=1}^M(\tilde{x}_t^m-\tilde{\mu}_t)(\tilde{x}_t^m-\tilde{\mu}_t)^T. 
\end{equation}

As a result, Eq.~(\ref{e:minkld}) is modified to %minimizing %$\kld(\pi_T,\hat{\pi}(x_t|y_{1:t}))$ where  is 
%\begin{subequations}
\begin{equation}
\min_{T\in\mathcal{H}} \kld(\pi_T,\hat{\pi}(x_t|y_{1:t})), \quad \mathrm{with}\quad
\hat{\pi}(\cdot|y_{1:t})\propto \hat{\pi}(\cdot|y_{1:t-1}) \pi(y_t|x_t).\label{e:minkld2}
\end{equation} 	
%\end{subequations}
Namely, we seek to minimize the distance between $\pi_T$ and the approximate posterior $\hat{\pi}(x_t|y_{1:t})$. 
We refer to the filtering algorithm by solving Eq.~\eqref{e:minkld2} as VEnKF,
where the complete algorithm is given in Alg.~\ref{alg:venkf}. 
\begin{algorithm}
	\caption{Affine-mapping based variational ensemble Kalman filter (AM-VEnKF)} \label{alg:venkf}
	\begin{itemize}
		\item Prediction: 
		\begin{itemize}
			\item Let $\tilde{x}_t^m\sim f_t(\cdot|x_{t-1}^m), m=1,2,\ldots,M$;
			\item Let $\hat{\pi}(\cdot|y_{1:t-1}) = N(\tilde{\mu}_t,\tilde{\Sigma}_t)$ where $\tilde{\mu}_t$ and $\hat{\Sigma}_t$
			are computed using Eq.~(\ref{eq:emp});
		\end{itemize}
		\item Update:
		\begin{itemize}
			\item Let $\hat{\pi}(x_t|y_{1:t})\propto \hat{\pi}(x_t|y_{1:t-1})\pi(y_t|x_t)$;
			\item Solve the minimization problem:
			\begin{equation}
			T_t=\arg\min_{T\in\mathcal{H}} \kld(\pi_T,\hat{\pi}(x_t|y_{1:t})).\nonumber
			\end{equation}
			\item Let
			$x_t^m=T_t\tilde{x}_t^m$ for $m=1,\ldots,M$.
		\end{itemize}
	\end{itemize}
\end{algorithm}

Now a key issue is to specify a suitable function space $\mathcal{H}$. %and loosely speaking, we choose $T$ to be an invertible affine mapping.
First let $A$ and $b$ be $n_x\times n_x$ and $n_x\times 1$ matrices respectively, 
and we can define a space of affine mappings
$\mathcal{A} = \{ T: T\cdot= A\cdot+b\}$,
with norm $\|T\| = \sqrt{\|A\|_2^2+\|b\|_2^2}$.
Now we choose 
\[ \mathcal{H} = \{ T \in \mathcal{A} \, |\, \|T\|\leq r, \,\mathrm{rank}(A) = n_x\},
\]
where $r$ is any fixed positive constant.
It is obvious that $A$ being full-rank implies that $T$ is invertible, which is an essential requirement for the proposed method,
and  will be discussed in detail in Section~\ref{sec:lof}.
%There are three main advantages for such a choice:
%\begin{itemize}
%\item It can be shown that a minimizer of Eq.~\eqref{e:minkld2} exists;
%\item It is easy to impose inevitability for the mappings;
%\item It can ensure that the resulting distributions do not deviate strongly from a Gaussian.
%\end{itemize}
Next we show that the minimizer of KLD exists in the closure of $\mathcal{H}$: 
\begin{thm}
	Let $P$ and $Q$ be two arbitrary probability distributions defined on $\mathcal{B}(\mathbb{R}^{n_x})$,
	%$\mathcal{A}$ be the space of all affine mappings from $\mathbb{R}^{n_x}$ to $\mathbb{R}^{n_x}$. 
	and 
	\begin{equation}
	{\mathcal{H}^*}=\{ T \in \mathcal{A} \, |\, \|T\|\leq r\},\nonumber
	\end{equation}
	for some fixed $r>0$.
	Let $P_T$ be the distribution of $T(x)$, given that $x$ be a $\mathbb{R}^{n_x}$-valued random variable following $P$. 
	The functional 
	$ \kld(P_T,Q)$  on $\mathcal{H}^*$ %\quad \mathrm{with}\quad
	admits a minimizer. %in $\mathcal{H}^*$. 
\end{thm}
\begin{proof}
	Let %$P_T=\nu(T)$ be a mapping from $\matchal{A}$ to 
	$\Omega=\{P_T:T\in\mathcal{H}^*\}$ be the image of $\mathcal{H}^*$ into $\mathcal{P}(\mathbb{R}^{n_x})$, the space of all Borel probability measures on $\mathbb{R}^{n_x}$. %and $ \nu(T) = P_T: \mathcal{H}^* \rightarrow \Omega$. 
	For any $\{T_n\} \in \mathcal{H}^*$ and $T\in \mathcal{H}^*$ such that
	$T_n\rightarrow T$, we have that
	$T_n(x)\rightarrow T(x)$ (a.s.), which implies that $P_{T_n}$ converges to $P_{T}$ weakly. It follows directly that $P_T$ is continuous on $\mathcal{H}^*$.
	%(Here $p_T:\mathcal{H}^*\rightarrow p(\mathbb{R}^{n_x})$, $p(\mathbb{R}^{n_x})$ endowed with the weak topology). 
	%  Since continuous mappings carry compact sets to compact sets, 
	Since $\mathcal{H}^*$ is a compact subset of $\mathcal{A}$, its image $\Omega$ is compact in $\mathcal{P}(\mathbb{R}^{n_x})$. 
	Since $\kld(P_T,Q)$ is lower semi-continuous with respect to $P_T$ (Theorem 1 in~\cite{posner1975random}), $\min\limits_{P_T\in\Omega} \kld(P_T,Q)$ 
	admits a solution $P_{T^*}$ with $T^*\in \mathcal{H}^*$.
	It follows that $T^*$ is a minimizer of $\min\limits_{T\in\mathcal{H}^*} \kld(P_T,Q)$.
\end{proof}
%\begin{cor} 
%Let $T^*$ be a solution of $\min_{T\in\mathcal{H}^*} \kld(p_T,q)$.
%There exists a minimizing sequence $\{T_n\} \in \mathcal{H}$ such that 
%$T_n \rightarrow T^*$ as $n\rightarrow \infty$. 
%%For any given $\epsilon>0$, there exist a mapping $T_\epsilon\in \mathcal{H}$ such that 
%%$T_\epsilon - T^*\| < \epsilon$. 
%\end{cor}
%\begin{proof}
%It follows directly from that $\mathcal{H}$ is dense in $\mathcal{H}^*$. 
%\end{proof}
%Loosely speaking, the states that, we are able to find an approximate solution in $\mathcal{H}$  that can be arbitrarily close to the actual minimizer within 
%its closure $\mathcal{H}^*$.
%from the computational perspective, since filtering often needs to be done sequentially and in realtime, the computational efficiency of a filtering algorithm is %essential. To this end, solving the optimization problem Eq.~(\ref{e:minkld}) for a general function space can pose a serious computational challenge, while %the use of affine mappings may considerably simplify the computation for solving Eq.~(\ref{e:minkld}). 

Finally it is also worth mentioning that, a key assumption of the proposed method (and EnKF as well) is that both the prior and posterior ensembles should not deviate strongly from Gaussian. To this end, a natural requirement for the chosen function space $\mathcal{H}$ is that, for any $T\in\mathcal{H}$, if $\pi(\tilde{x}_t|y_{1:t-1})$ is close to Gaussian, so should be $\pi_T(x_t)$ with $x_t=T(\tilde{x}_t)$. 
Obviously an arbitrarily function space does not satisfy such a requirement. However, for affine mappings, we have the following proposition: 
\begin{prop}
	For a given positive constant number $\epsilon$,
	if there is a $n_x$-dimensional normal distribution $\tilde{p}_G$ such that 
	$\kld(\tilde{p}_G(\tilde{x}_t),\pi(\tilde{x}_t|y_{1:t-1}))<\epsilon$, and if $T \in \mathcal{H}$, there must exist a $n_x$-dimensional normal distribution ${p}_G$ satisfying $\kld({p}_G({x}_t),\pi_T(x_t))<\epsilon$. 
\end{prop}
\begin{proof}
	This proposition is a direct consequence of the fact that KLD is invariant under affine transformations.
\end{proof}

\noindent Loosely the proposition states that, for an affine mapping $T$, if the prior $\pi(\tilde{x}_t|y_{1:t-1})$ is close to a Gaussian distribution, so is $\pi_T(x_t)$, which ensures that the update step will not increase the ``non-Gaussianity'' of the ensemble.

{%\color{red}
	%Another key advantage of  the affine mappings is that one can rather easily impose the invertibility condition on them, which is an essential theoretical %requirement  of the method; this will be discussed in detail in Section~\ref{sec:lof}. 
	In principle one can choose a different function space $\mathcal{H}$, and for example, a popular transport-based approach called the Stein variational gradient descent (SVGD) method \cite{liu2016stein} constructs such a function space using the reproducing kernel Hilbert space (RKHS), which can also be used in the VEnKF formulation. We provide a detailed description of the SVGD based VEnKF in Appendix~A, and this method is also compared with the proposed AM-VEnKF in all the numerical examples.}

%	Finally we consider a special case where distribution $\pi(x_{t-1}|y_{1:t-1})$ from the previous step is Gaussian,
%	and the propagation model is in the form of
%	$x_{t} = F_t x_{t-1} +\eta_t,$ 
%	where $F_t$ is a $d\times d$ matrix and $\eta_t$ follows a zero-mean Gaussian distribution with variance $\Sigma^\eta_t$. In this case we do not need to represent the posterior distributions with ensembles,
%	and instead we can derive a affine mapping based Kalman filter which computes the posteriors directly. 

%	\begin{algorithm}
%	\caption{The Linear-mapping based Kalman filter (LMKF)}
%		\begin{itemize}
%		\item Prediction: 
%		\begin{itemize}
%		\item Let $\tilde{\mu}_t= F_t\mu_{t-1}$ and $\tilde{\Sigma}_t =F_t\Sigma_{t-1}F_t+\Sigma_t^\eta$;
%		\item Let $\pi(\tilde{x}_t|y_{1:t-1})=N(\tilde{\mu},\tilde{\Sigma}_t)$;
%		\end{itemize}
%		      \item Update:
%		      
%		      \begin{itemize}
%		         \item Let ${\pi}(x_t|y_{1:t})\propto {\pi}(x_t|y_{1:t-1})g_t(y_t|x_t)$;
%
%                \item Solve the minimization problem:
%		          
%        	\begin{equation}
%	       T_t=\arg\min_{T\in\mathcal{L}} \kld(\pi_T,{\pi}(x_t|y_{1:t}));
%	\end{equation} 
%		\item Let $\mu_t=T_t\tilde{\mu}$ and $\Sigma_t=T_t\tilde{\Sigma}_tT_t^{T}$.
%	\end{itemize}
%	\end{itemize}
%	\end{algorithm}
%	

\subsection{Connection to the ensemble Kalman filter} \label{sec:enkf}
In this section, we discuss the connection between the standard EnKF and AM-VEnKF,
and show that EnKF results in additional estimation error due to certain approximations made. 
We start with a brief introduction to  EnKF. 
We consider the situation where the observation model takes the form of
% \[g_t(y_t|x_t)=N(H_t x_t,R_t).\]
{
	\begin{equation}\label{eq:lom}
	y_t=\textsc{H}_tx_t+\beta_t,
	\end{equation}
	which implies $\pi(y_t|x_t)=N(H_t x_t,R_t)$,
	where $H_t$ is a linear observation operator and $\beta_t$ is a zero-mean Gaussian noise with covariance $R_t$. } 

In this case, EnKF can be understood as to obtain an approximate solution of Eq.~\eqref{e:minkld2}. 
Recall that in the VEnKF formulation,  $\pi_T$ is the distribution of 
$x_t = T(\tilde{x}_t)$ 
where  $\tilde{x}_t$ follows $\pi(\cdot|y_{1:t-1})$,
and similarly we can define  $\hat{\pi}_T$ as the distribution of 
$x_t = T(\tilde{x}_t)$ where  $\tilde{x}_t$ follows the approximate prior $\hat{\pi}(\cdot|y_{1:t-1})$.
Now instead of Eq.~\eqref{e:minkld2}, we find $T$ by solving,
\begin{equation}
\min_{T\in\mathcal{H}} \kld(\hat{\pi}_T,\hat{\pi}(x_t|y_{1:t})),\label{e:minkld3}
\end{equation} 	
and the obtained mapping $T$ is then used to transform the particles. 
It is easy to verify that the optimal solution of Eq.~(\ref{e:minkld3}) can be obtained exactly, 
\begin{equation}
x_t = T(\tilde{x}_t)= (\mathrm{I}-K_t H_t)\tilde{x}_t+K_ty_t,\label{e:enkfmap}
\end{equation}
where $\mathrm{I}$ is the identity matrix and Kalman Gain matrix $K_t$ is
\begin{equation}
\label{eq:kgain}
K_t=\tilde{\Sigma}_tH_t^T(H_t\tilde{\Sigma}_tH_t^T+\mathrm{R}_t)^{-1}.
\end{equation}
Moreover, the resulting value of KLD is zero, which means that the optimal mapping pushes the prior  exactly to the posterior. One sees immediately that the optimal mapping in Eq.~(\ref{e:enkfmap}) coincides with the updating formula of EnKF, 
implying that EnKF is an approximation of VEnKF, even when the observation model is exactly linear-Gaussian. 
%Next we show that this approximation is only sub-optimal. 

%\begin{thm}
%Let $T_{EnKF}$ be the mapping given by Eq.~\eqref{e:enkfmap}. There exists an invertible linear mapping $T$ such that
%\[ D_{KL} \pi_T \leq D_{KL}, \]
%where the equality sign holds iff $\pi=\hat{\pi}$. 
%\end{thm}

%it has been assumed in LM-VEnKF that the prior is approximated as $\hat{\pi}(\cdot|y_{1:t-1})=N(\tilde{\mu}_t,\tilde{\Sigma}_t)$, and it follows that the %approximate posterior is also Gaussian: $\hat{\pi}(\cdot|y_{1:t}) = N(\mu_t,\Sigma_t).$
%The mean $\mu_t$ and the covariance $\Sigma_t$ can be obtained analytically: 
%	\begin{equation}
%	\mu_t=(\mathrm{I}-K_tH_t)\tilde{\mu}_t+K_ty_t,\quad  \Sigma_t=(\mathrm{I}-K_tH_t)\tilde{\Sigma}_t.
%	\end{equation}
%	where $\mathrm{I}$ is the identity matrix and Kalman Gain matrix $K_t$ is
%	\begin{equation}
%	\label{eq:kgain}
%	K_t=\tilde{\Sigma}_tH_t^T(H_t\tilde{\Sigma}_tH_t^T+\mathrm{R}_t)^{-1}.
%	\end{equation}

%Following similar steps, one can also derive the standard Kalman filter as a special case of LMKF for a linear-Gaussian observation model.
%Finally we note that a simple extension of EKF to the generic observation model is available where the Kalman gain matrix is obtained 
%via a sample estimate of the covariance between the state vector and the observation vector~\cite{houtekamer2001sequential}. This extended EKF method will be compared against the proposed LMEKF in the numerical experiments.

When the observation model is not linear-Gaussian, further approximation is needed. 
Specifically the main idea  is to approximate the actual observation model with a linear-Gaussian one, and estimate the Kalman gain matrix $K_t$ directly from the ensemble~\cite{houtekamer2001sequential}. Namely, suppose we have an ensemble from the prior distribution: $\{\tilde{x}_t^m\}_{m=1}^M$, and we generate an ensemble of data points: $\tilde{y}_t^m\sim \pi(\tilde{y}_t^m|\tilde{x}_t^m)$ for $m=1,\ldots,M$. Next we estimate the Kalman gain matrix as follows,
\begin{eqnarray*}
	&\tilde{K}_t=C_{xy}C_{yy}^{-1},\\
	&\hat{x} _t= \frac1M\sum_{m=1}^M \tilde{x}^m_t,\quad\hat{y}_t = \frac1M\sum_{m=1}^M \tilde{y}^m_t,\\
	&C_{xy}=\frac{1}{M-1}\sum\limits_{m=1}^{M}(\tilde{x}_t^m-\hat{x}_t)(\tilde{y}_t^m-\hat{y}_t)^T,\\
	&C_{yy}=\frac{1}{M-1}\sum\limits_{m=1}^M(\tilde{y}_t^m-\hat{y}_t)(\tilde{y}_t^m-\hat{y}_t)^T.
\end{eqnarray*}

Finally the ensemble are updated: $x_t^m=\tilde{x}_t^m+\tilde{K}_t(y_t-\tilde{y}_t^m)$ for $i=1,\ldots,M$. 
As one can see here, due to these approximations, the EnKF method can not provide an accurate solution
to Eq.~\eqref{e:minkld2}, especially when these approximations are not  accurate. 
%This version of EnKF will be compared against the proposed AM-VEnKF in the numerical experiments. 

\subsection{Numerical algorithm for minimizing KLD}
\label{sec:lof}
In the VEnKF framework presented in section~\ref{sec:lmkf}, the key step is to solve KLD minimization problem~(\ref{e:minkld}). In this section we describe in details how the optimization problem is solved numerically. 

Namely suppose at step $t$, we have a set of samples $\{\tilde{x}^m_{t}\}_{m=1}^M$ drawn from the prior distribution $\pi(\tilde{x}_t|y_{1:t-1})$, we want to transform them into the ensemble $\{{x}^m_{t}\}_{m=1}^M$ that follows the approximate posterior $\pi({x}_t|y_{1:t})$. First we set up some notations, and for conciseness some of them are different from those used in the previous sections: first we drop the subscript of $\tilde{x}_t$ and $x_t$, and we then define ${p}(\tilde{x})={\pi}(\tilde{x}|y_{1:t-1})$ (the actual prior), %$\hat{p}(\tilde{x})=\tilde{\pi}(\tilde{x}|y_{1:t-1})=N(\tilde{\mu},\tilde{\Sigma})$ (the Gaussian approximate prior), 
$\tilde{p}(\tilde{x})=\hat{\pi}(\tilde{x}|y_{1:t-1})=N(\tilde{\mu},\tilde{\Sigma})$ (the Gaussian approximate prior), $l(x)=-\log\pi(y_{t}|x)$ (the negative log-likelihood) and 
$q(x) = \hat{\pi}(x|y_{1:t})$ (the approximate posterior).
It should be clear that
\begin{equation}
q(x) \propto \tilde{p}(\-x) \exp(-l(x)). \label{e:apppost}
\end{equation}

Recall that we want to minimize $\kld(p_T(x),q(x))$ where $p_T$ is the distribution of the transformed random variable $x=T(\tilde{x})$, and it is easy to show that 
\begin{equation}
\kld(p_T(x),q(x)) = \kld(p(\tilde{x}),q_{T^{-1}}(\tilde{x})),\nonumber
\end{equation}
where $q_{T^{-1}}$ is the distribution of the inversely transformed random variable $\tilde{x}=T^{-1}(x)$ with $x\sim q(x)$.
Moreover, as
\begin{equation}
\kld(p(\tilde{x}),q_{T^{-1}}(\tilde{x})) = \int \log[p(\tilde{x})]p(\tilde{x}) d\tilde{x} - \int \log[q_{T^{-1}}(\tilde{x})]p(\tilde{x}) d\tilde{x},\nonumber
\end{equation}
minimizing $\kld(p_T(x),q(x))$  is equivalent to 
\begin{equation}
\min_{T\in\mathcal{H}}-\int \log[q_{T^{-1}}(\tilde{x})]p(\tilde{x}) d\tilde{x}. \label{e:max1}
\end{equation}
A difficulty here is that the feasible space $\mathcal{H}$ is constrained by $\|T\|\leq r$ (i.e. an Ivanov regularization), which poses computational challenges.
%Eq.~\eqref{e:max1} as a constrained optimization problem.
Following the convention we replace the constraint with a Tikhonov regularization to simplify the computation:
\begin{equation}
\min_{T\in\mathcal{A}}-\int \log[q_{T^{-1}}(\tilde{x})]p(\tilde{x}) d\tilde{x}+\lambda \|T\|^2, \label{e:min_reg}
\end{equation}
where $\lambda$ is a pre-determined regularization constant.
%Since the function space $\mathcal{L}$ represents affine and bijective mappings, we can write it as
%\[\mathcal{L}=\{Tx = Ax+b\,|\, \mathrm{rank}(A)=n_x,\,b\in \mathbb{R}^{n_x}\},\]
%and thus we just need determine the matrix $A$ and the vector $b$. 

%Next we substitute $Tx=Ax+b$ into Eq.~\eqref{e:min_reg}, yielding, 
Now using $Tx=Ax+b$, $q_{T^{-1}}(\tilde{x})$ can be written as, 
\begin{equation}
q_{T^{-1}}(\tilde{x})=q(A\tilde{x}+b)|A|, \label{e:qT-1}
\end{equation}
and we substitute Eq.~(\ref{e:qT-1}) along with Eq.~(\ref{e:apppost}) in to Eq.~(\ref{e:min_reg}), yielding,
\begin{eqnarray}
\min_{A,b}F_q(A,b)%&&\nonumber\\
&:=&-\int \log[q(A\tilde{x}+b)]p(\tilde{x}) d\tilde{x} -\log|A|+\lambda(\|A\|_2^2+\|b\|_2^2), \nonumber \\
&=&-\int \log [\tilde{p}(A\tilde{x}+b)] p(\tilde{x})d\tilde{x} + \int l(A\tilde{x}+b) p(\tilde{x}) d\tilde{x}\nonumber\\
&&\quad -\log|A|+\lambda(\|A\|_2^2+\|b\|_2^2),\nonumber\\
&=&\frac{1}{2}Tr[(\tilde{\Sigma}+\tilde{\mu}\tilde{\mu}^{T})A^T\tilde{\Sigma}^{-1}A]
+(b-\tilde{\mu})^T{\tilde{\Sigma}^{-1}}[A\tilde{\mu}+\frac{1}{2}(b-\tilde{\mu})]\nonumber\\
&&\quad-\log|A|+\mathrm{E}_{\tilde{x}\sim p}[l(A\tilde{x}+b)]+\frac{1}{2}(n_x\log(2\pi)+\log{|\tilde{\Sigma}|})\nonumber\\
&&\quad +\lambda(\|A\|_2^2+\|b\|_2^2),\label{e:maxF}
\end{eqnarray}
which is an unconstrained optimization problem in terms of $A$ and $b$. 
It should be clear that the solution of Eq.~\eqref{e:maxF} is naturally invertible. 

We then solve the optimization problem~(\ref{e:maxF}) with a gradient descent  (GD) scheme: 
\begin{eqnarray*}
	A_{k+1} &=& A_k-\epsilon_k   \frac{\partial F_q}{\partial A}(A_k,b_k),\\
	b_{k+1} &=& b_k-\epsilon_k  \frac{\partial F_q}{\partial b}(A_k,b_k),
\end{eqnarray*}
where $\epsilon_k$ is the step size and the gradients can be derived as,
%\numparts
\begin{eqnarray}\label{e:grads}
\frac{\partial F_q}{\partial A}(A,b)&=& (\tilde{\Sigma}+\tilde{\mu}\tilde{\mu}^{T})A^T{\tilde{\Sigma}^{-1}}+{\tilde{\Sigma}^{-1}}(b-\tilde{\mu})\tilde{\mu}^T -A^{-1} 
\nonumber
\\&& +\mathrm{E}_{\tilde{x}\sim p}[ \nabla_xl(A\tilde{x}+b)\tilde{x}^T]+2\lambda A, \\
\frac{\partial F_q}{\partial b}(A,b)&=&\tilde{\Sigma}^{-1}[A\tilde{\mu}+b-\tilde{\mu}]+\mathrm{E}_{\tilde{x}\sim p}[ \nabla_xl(A\tilde{x}+b)]
+2\lambda b.
\end{eqnarray}
%\endnumparts

Note that Eq.~(\ref{e:grads}) involves the expectations $\mathrm{E}_{\tilde{x}\sim p}[ \nabla_xl(A\tilde{x}+b)\tilde{x}^T]$ and $\mathrm{E}_{\tilde{x}\sim p}[\nabla_xl(A\tilde{x}+b)]$ which are not known exactly, and in practice they can be replaced by their Monte Carlo estimates:
\begin{eqnarray*}
	&\mathrm{E}_{\tilde{x}\sim p}[ \nabla_xl(A\tilde{x}+b)\tilde{x}^T] \approx\frac1M \sum   \nabla_xl(A\tilde{x}^m+b)(\tilde{x}^m)^T,\\
	&\mathrm{E}_{\tilde{x}\sim p}[ \nabla_xl(A\tilde{x}+b)]\approx 
	\frac1M\sum_{m=1}^M \nabla_xl(A\tilde{x}^m+b),
\end{eqnarray*}
where $\{\tilde{x}^m\}_{m=1}^M$ are the prior ensemble and $\nabla_xl(x)$ is the derivative  of $l(x)$ taken with respect to $x$. 
The same Monte Carlo treatment also applies to the objective function $F_q(A,b)$ itself when it needs to be evaluated. 
%{\color{red} It is worth noting that, for some specific likelihood functions, the optimization can be simplified as the resulting objective function may be of some special structure.
%Since we consider the generic models in this work, such an option is not explored here.}

The last key ingredient of the optimization algorithm is the stopping criteria. Due to the stochastic nature of the optimization problem, standard stopping criteria in the gradient descent method are not effective here. Therefore we adopt a commonly used criterion in search-based optimization: the iteration is terminated if the current best value is not sufficiently increased within a given number of steps. More precisely, let $F^*_k$ and $F^*_{k-\Delta k}$ be the current best value at iteration $k$ and $k-\Delta k$ respectively where $\Delta k$ is a positive integer smaller than $k$, and the iteration is terminated if $F^*_{k}-F^*_{k-\Delta k}< \Delta_F$ for a prescribed threshold $\Delta_F$. In addition we also employ a safeguard stopping condition, which terminates the procedure after the number of iterations reaches a prescribed value $K_{\max}$. 

It is also worth mentioning  that the EnKF type of methods are often applied to problems where the ensemble size is similar to or even smaller than the dimensionality of the states and in this case the  localization techniques are usually used to address the undersampling issue~\cite{anderson2007exploring}. In the AM-VEnKF method, many localization techniques developed in EnKF literature can be directly used, and in our numerical experiments we adopt the sliding-window localization used in~\cite{ott2004local}, and we will provide more details of this localization technique in Section~\ref{sec:obsmodel}. 
%where local observations are used to update local state vectors, and the whole state vector is reconstructed by aggregating the local updates.

%
%
%	\begin{algorithm}
%		\caption{Gradient Descent} %
%		\label{alg:GD}
%		\begin{algorithmic}[1]
%			\STATE{\bf{Input:} A target function $F_q(A,b)$}, initial particles $\{x_t^m\}_{m=1}^M$, $A_0=I,b_0=0$,and threshold value $s$.
%			\STATE{\bf{Output:} A set of updated particles $\{x_{t,I}^m\}_{m=1}^M$ and the iterations $I$.}
%			\FOR{ iteration $i=1$ }
%			\STATE{$A_i=A_{i-1}+\epsilon_i\frac{\partial F_q}{\partial A}$,\\
%				$b_i=b_{i-1}+\epsilon_i\frac{\partial F_q}{\partial b}$,\\ 
%				and $\epsilon_i$ is the stepsize at the $i$-th iteration.}
%			\IF{$||x_{t,i}-x_{t,i-1}||<s$}
%			\STATE{Stop the iteration, record $I=i$, return particles $\{x_{t,I}^m\}_{m=1}^M$;}
%			\ELSE
%			\STATE{i=i+1;}
%			\ENDIF
%			\ENDFOR
%		\end{algorithmic}
%	\end{algorithm}

\section{Numerical examples}\label{sec:examples}
\subsection{Observation models}\label{sec:obsmodel}

As is mentioned earlier, the goal of this work is to deal with generic observation models, and in our numerical experiments, we test the proposed method with { an observation model that is quite flexible and also commonly used in epidemic modeling and simulation~\cite{capaldi2012parameter}:
	\begin{equation}\label{eq:nn}
	y_t=G(x_t,\beta_t)=M(x_t)+aM(x_t)^{\theta}\circ\beta_t,
	\end{equation}
	where $M(\cdot): \mathcal{X}\rightarrow \mathcal{Y}$ is a mapping from the state space to the observation space}, $a$ is a positive scalar, $\beta_t$ is a random variable defined on $\mathcal{Y}$, and $\circ$ stands for the Schur (component-wise) product. Moreover we assume that $\beta_t$ is an independent random variable with zero mean and variance $R$, where $R$ here is the vector containing the variance of each component and should not be confused with the covariance matrix. It can be seen that $a M(x_t)^{\theta}\circ\beta_t$ represents the observation noise, controlled by  two adjustable parameters $\theta$ and $a$, and the likelihood $\pi(y_t|x_t)$ is of mean $M(x_t)$ and variance $a^2M(x_t)^{2\theta}\circ R$.

The parameter $\theta$ is particularly important for specifying the noise model in~\cite{capaldi2012parameter} and here we consider the following three representative cases. First if we take $\theta=0$, it follows that $y_t=M(x_t)+a\beta_t$, where the observation noise is  independent of the state value $x_t$. This is the most commonly used observation model in data assimilation and we refer to it as the absolute noise following~\cite{capaldi2012parameter}. Second if $\theta=0.5$, the variance of observation noise is $a^2M(x_t)\circ R$, which is linearly dependent on $M(x_t)$, and we refer to this as the Poisson noise~\cite{capaldi2012parameter}. Finally in case of $\theta=1$, it is the standard deviation of the noise, equal to $aM(x_t)R^{1/2}$, that depends linearly on $M(x_t)$, and this case  is referred to as the relative noise~\cite{capaldi2012parameter}. In our numerical experiments we test all the three cases. 

Moreover, in the first two numerical examples provided in this work, we take 
\begin{equation}
M(x_t) =0.1x^2_t,\label{e:M}
\end{equation}
$a=1$,
and assume $\beta_t$ to follow the Student's $t$-distribution~\cite{6638770} with zero-mean and variance 1.5. In the last example, we take,
\begin{equation}
M(x_t) =\exp(x_t/2),\label{e:M2}
\end{equation}
and $a=1$.

As has been mentioned, localization is needed in some numerical experiments here. Given Eqs.~(\ref{e:M}) and (\ref{e:M2}) we can see that the resulting observation model has a property that each component of the observation $y_t$ is associated to a component of the state $x_t$: namely,
\[y_{t,i} = M(x_{t,i})+(M(x_{t,i}))^\theta\beta_{t,i},\quad i=1,\ldots,n_x,\]
where $\beta_{t,i}$ is the $i$-th component of $\beta_t$, and $n_y=n_x$. In this case, we can employ the sliding-window localization method, where local observations are used to update local state vectors, and the whole state vector is reconstructed by aggregating the local updates. Namely, the state vector $x_t=(x_{t,1},\ldots,x_{t,n_x})$ is decomposed into a number of overlapping local vectors: $\{x_{t,N_i}\}_{i=1}^{n_x}$, 
% $$\Omega=\{x_{t,1:l+i}\}_{i=1}^l\cup\{x_{t,i-l:i+l}\}_{i=l+1}^{d-l}\cup\{x_{t,i-l:d}\}_{i=d-l+1}^d,$$
where $N_i = [\max\{1,i-l\}: \min\{i+l,n_x\}]$
% \[ N_i = \left\{\begin{matrix}
% 1:l+i&i\\
% i-l:i:l&i\\
% i-l:d,& d-l+1\leq i\leq d
% \end{matrix}\right.
% \]
for a positive integer $l$.
%For simplicity, we define the vectors in $\Omega$ as $x_{t,N_1}$,
When updating any local vector $x_{t,N_i}$, we only use the local observations $y_{t,N_i}$ and as such each local vector  is updated independently. It can be seen that by design each $x_{t,i}$ is updated in multiple local vectors, and the final update is calculated by averaging its updates in local vectors indexed by $N_{\max\{1,i-k\}},\ldots,N_{i},\ldots, N_{\min\{i+k,n_x\}}$, for some positive integer $k\leq l$. We refer to \cite{ott2004local,lei2011moment} for further details. 

\subsection{Lorenz-96 system}\label{sec:lor96}
Our first example is the Lorenz-96 model~\cite{lorenz1996predictability}:
\begin{equation} \label{eq:lorenz}
\left\{\begin{array}{ll}
\frac{dx^n}{dt}=(x^{n+1}-x^{n-2})x^{n-1}-x^{n}+8,\ n=1,\ldots,40\\
x^0=x^{40},\ x^{-1}=x^{39},\ x^{41}=x^1,
\end{array}
\right.
\end{equation}
a commonly used benchmark example for filtering algorithms.

By integrating the system~(\ref{eq:lorenz}) via the Runge-Kutta scheme with stepsize $\Delta t=0.05$, and adding some model noise, we obtain the following  discrete-time model:
\begin{equation} \label{eq:discretelorenz}
\left\{
\begin{array}{ll}
\textbf{x}_t &= \mathcal{F}(\textbf{x}_{t-1})+\alpha_t,\quad t=1,2,\ldots\\
\textbf{y}_t &= M(\textbf{x}_t)+M(\textbf{x}_t)^{\theta}\beta_t, \quad t=1,2,\ldots
\end{array}
\right.
\end{equation}
where $\mathcal{F}$ is the standard fourth-order Runge-Kutta solution of Eq.~(\ref{eq:lorenz}), $\alpha_t$ is standard Gaussian noise, and the initial state $\textbf{x}_0\sim U[0,10]$. We use synthetic data in this example, which means that both the true states and the observed data are simulated from the model. 

%I, $\beta_t$ is Student's t distribution with parameters $v$, while the t distribution reduces to the Gaussian as $v$ tends to infinity, thus includes it as a %special case, here we denote that $v=[6,\ldots,6]_{40\times1}$. And in this example, the stepsize in Gradient Descent scheme is 0.001 at every iteration. 

As mentioned earlier, we consider the three observation models corresponding to $\theta=0, 0.5$ and $1$. In each case, we use two sample sizes $M=100$ and $M=20$. To evaluate the performance of VEnKF, we implement both the AM based and the SVGD based VEnKF algorithms. As a comparison, we also impliment several commonly used methods: the EnKF variant provided in Section~\ref{sec:enkf}, PF, and NLEAF~\cite{lei2011moment} with first-order  (denoted as NLEAF 1) and second-order (denoted as NLEAF 2) correction, in the numerical tests. The stopping criterion in AM-VEnKF is specified by $\Delta_k=20$, $\Delta_F=0.1$ and $K_{\max}=1000$, while the step size $\epsilon_k$ in GD iteration is $0.001$. In SVGD-VEnKF, the step size is also  $0.001$, and the stopping criterion is chosen in a way so that the number of iterations is approximately the same as that in AM-VEnKF. For the small sample size $M=20$, in all the methods except PF, the sliding window localization (with $l=3$ and $k=2$; see \cite{lei2011moment} for details) is used. 

With each method, we compute the estimator bias (i.e., the difference between the ensemble mean and the ground truth) at each time step and then average the bias over the 40 different dimensions. The procedure is repeated 200 times for each method and all the results are averaged over the 200 trials to alleviate the statistical error. 

\begin{figure}[htbp]
	%\begin{minipage}[htbp]{0.5\linewidth}
	\centerline{\includegraphics[width=0.6\textwidth]{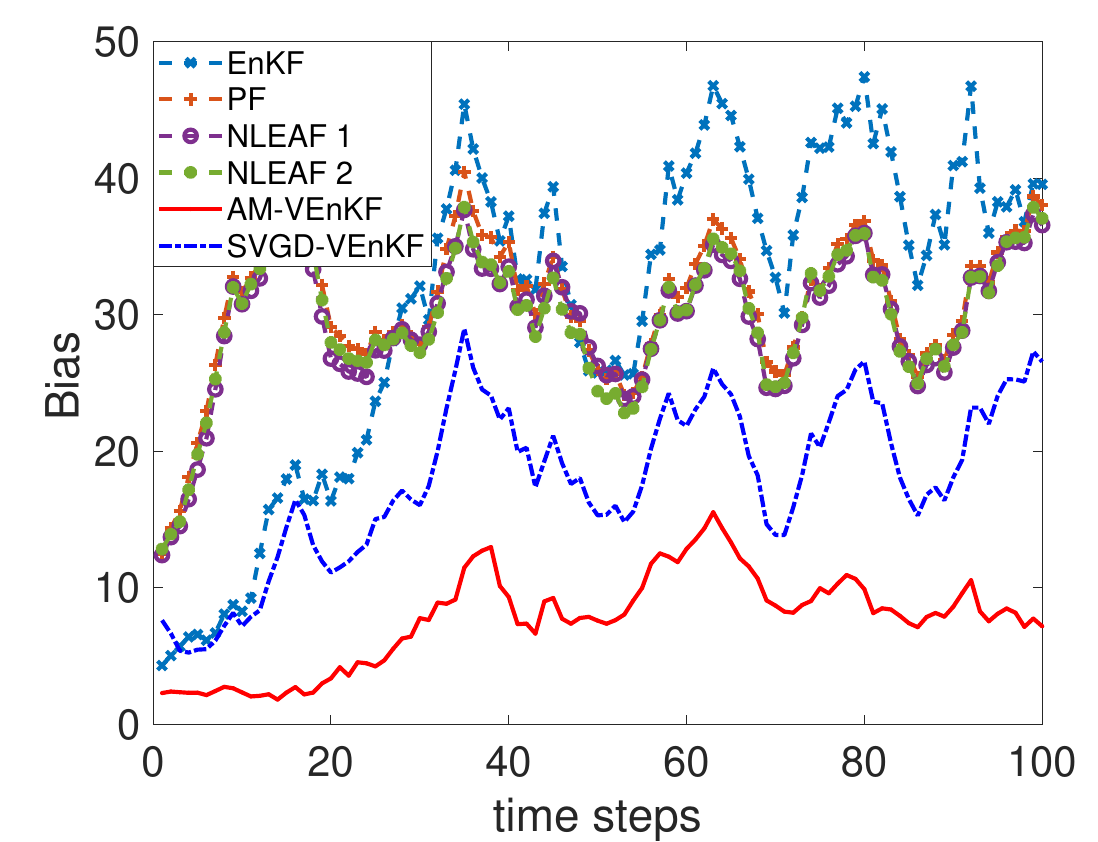}}
	\caption{The average bias at each time step for $\theta=0$ and $M=100$ in the Lorenz 96 example. }
	\label{fig:t1_b0}
	
	\centerline{\includegraphics[width=.5\textwidth]{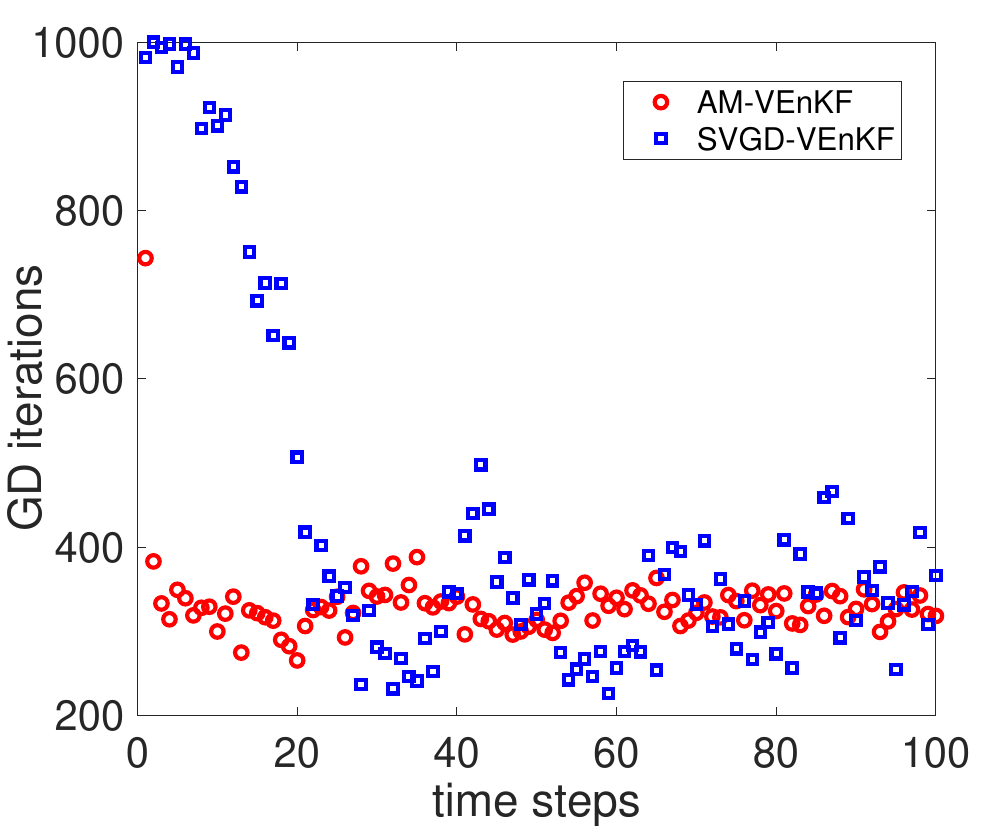}\includegraphics[width=.5\textwidth]{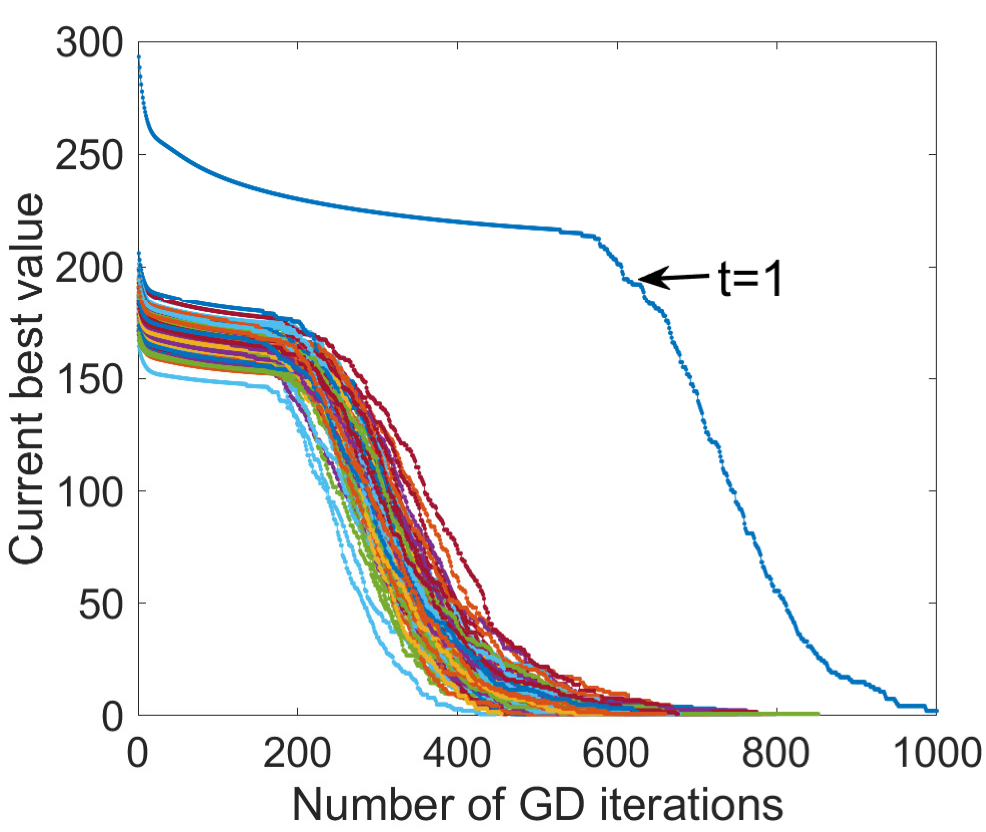}}
	\caption{Left: the number of GD iterations (in both AM and SVGD) at each time step. 
		Right: the current best value plotted against the GD iterations (in AM) where 
		each line represents a time step. The results are for $\theta =0$ and $M=100$ in the Lorenz 96 example.}
	\label{fig:t1_g0}
\end{figure}

The average bias for $\theta =0$ is shown in Fig.~\ref{fig:t1_b0} where it can be observed that in this case, while the other three methods yield largely comparable accuracy in terms of estimation bias, the bias of AM-VEnKF is significantly smaller. To analyze the convergence property of the method, in Fig.~\ref{fig:t1_g0} (left) we show the number of GD iterations (of both AM and SVGD) at each time step, where one can see that all GD iterations terminate after around 300-400 steps in AM-VEnKF, except the iteration at $t=1$ which proceeds for around 750 steps. The SVGD-VEnKF undergoes a much higher number of iterations in the first 20 time steps, while becoming about the same level as that of AM-VEnKF. This can be further understood by observing Fig.~\ref{fig:t1_g0} (right) which shows the current best value $F^*_k$ with respect to the GD iteration in AM-VEnKF, and each curve in the figure represents the result at a time step $t$. We see here that the current best values  become settled after around 400 iterations at all time locations except  $t=1$, which agrees well with the number of iterations shown on the left. It is sensible that the GD algorithm takes substantially more iterations to converge at $t=1$, as the posterior at $t=1$ is typically much far away from the prior, compared to other time steps. These two figures thus show that the proposed stopping criteria are effective in this example. 
\begin{figure}[htbp]
	\centering
	\includegraphics[width=0.6\textwidth]{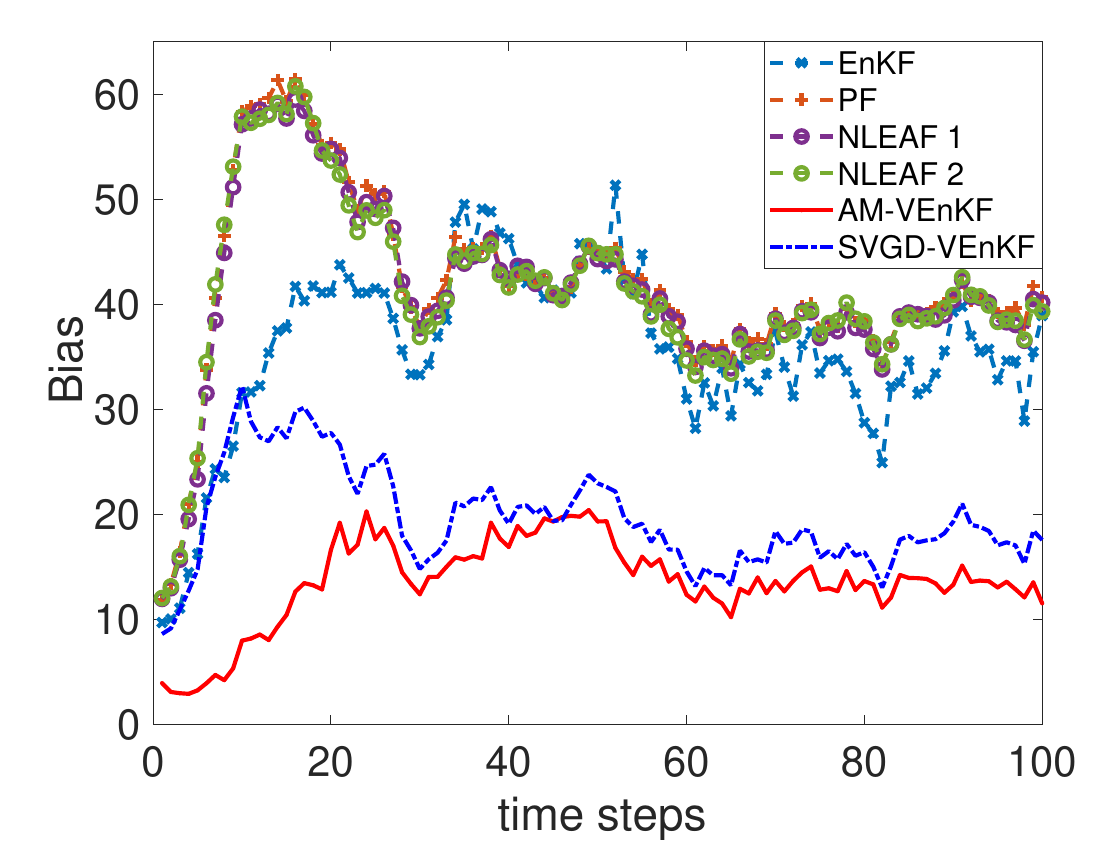}
	\caption{The average bias at each time step for $\theta=0.5$ and $M=100$ in the Lorenz 96 example. }
	\label{fig:t1_b05}
	\centerline{\includegraphics[width=.5\textwidth]{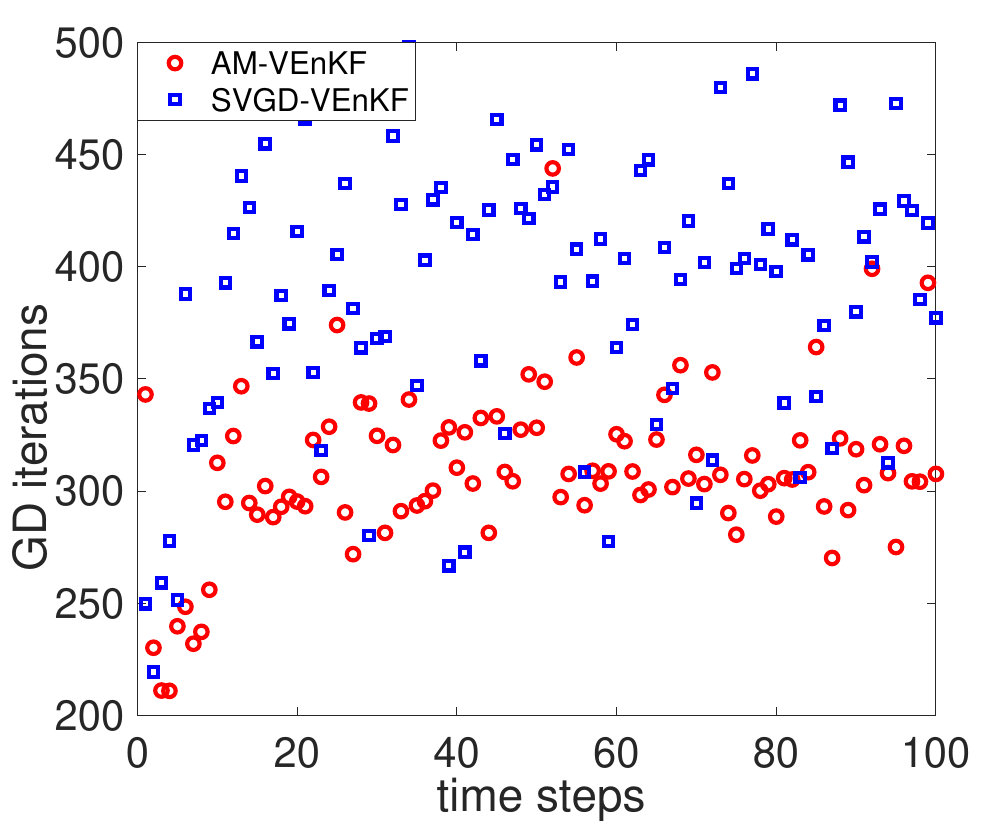}
		\includegraphics[width=.5\textwidth]{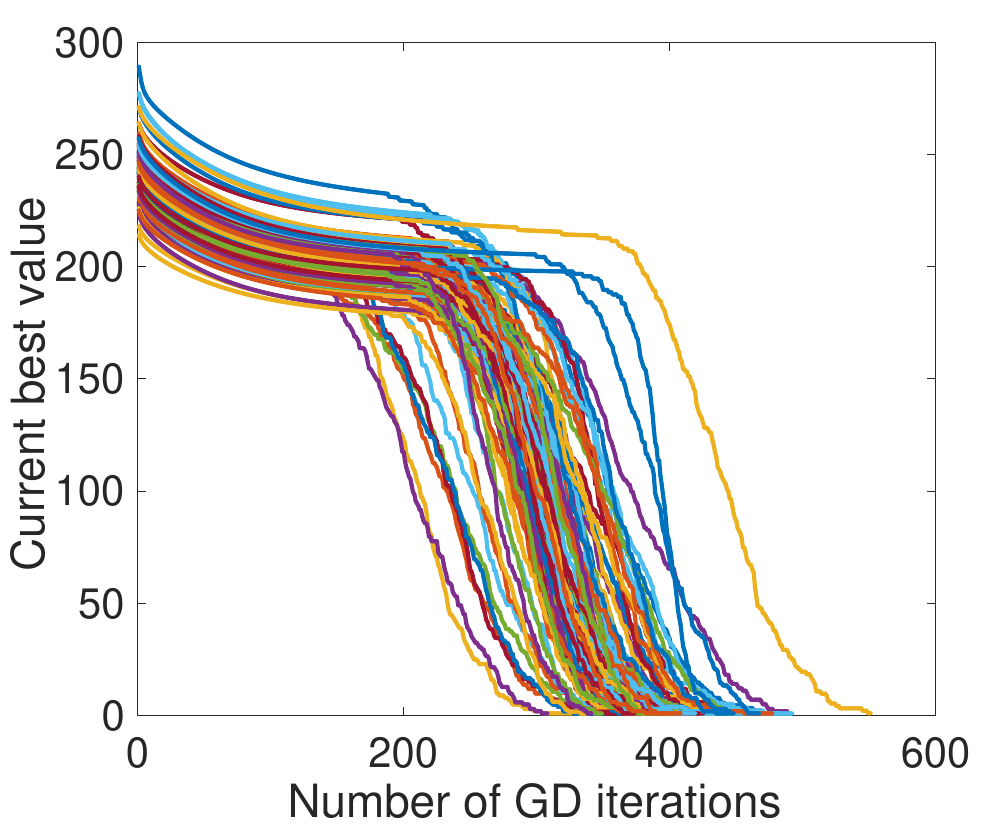}}
	\caption{Left: the number of GD iterations  (in both AM and SVGD) at each time step. 
		Right: the current best value plotted against the GD iterations (in AM) where 
		each line represents a time step. The results are for $\theta =0.5$ and $M=100$ in the Lorenz 96 example.}
	\label{fig:t1_g05}
\end{figure}

\begin{figure}[htbp]
	%\begin{minipage}[htbp]{0.5\linewidth}
	\centering
	\includegraphics[width=0.6\textwidth]{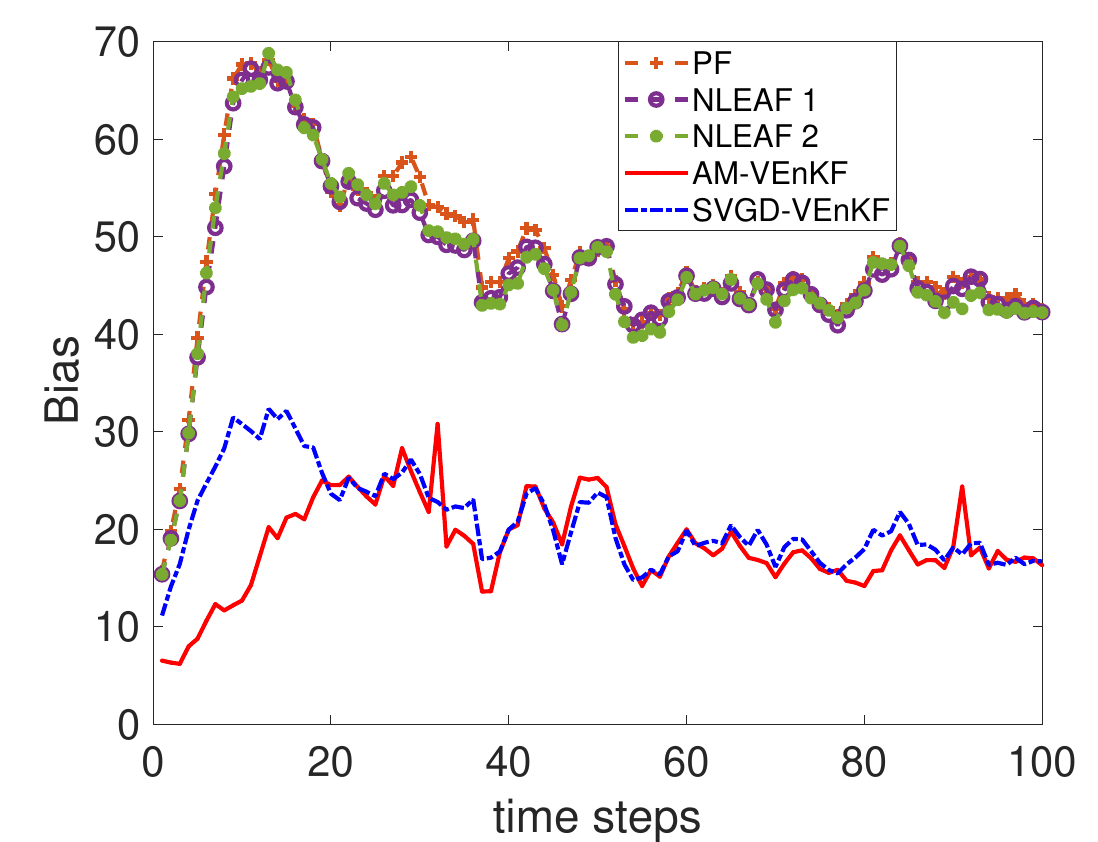}
	\caption{The average bias  at each time step for $\theta =1$ and $M=100$ in the Lorenz 96 example. }
	\label{fig:t1_b1}
	\centerline{\includegraphics[width=0.5\textwidth]{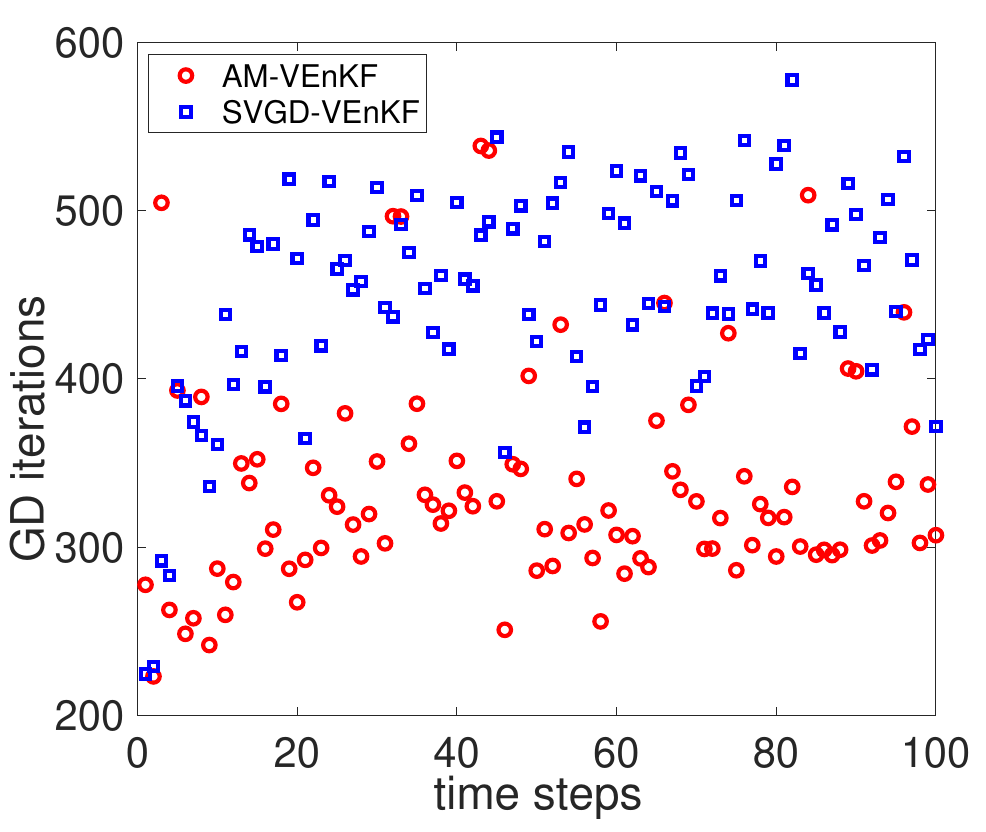}	
		\includegraphics[width=0.5\textwidth]{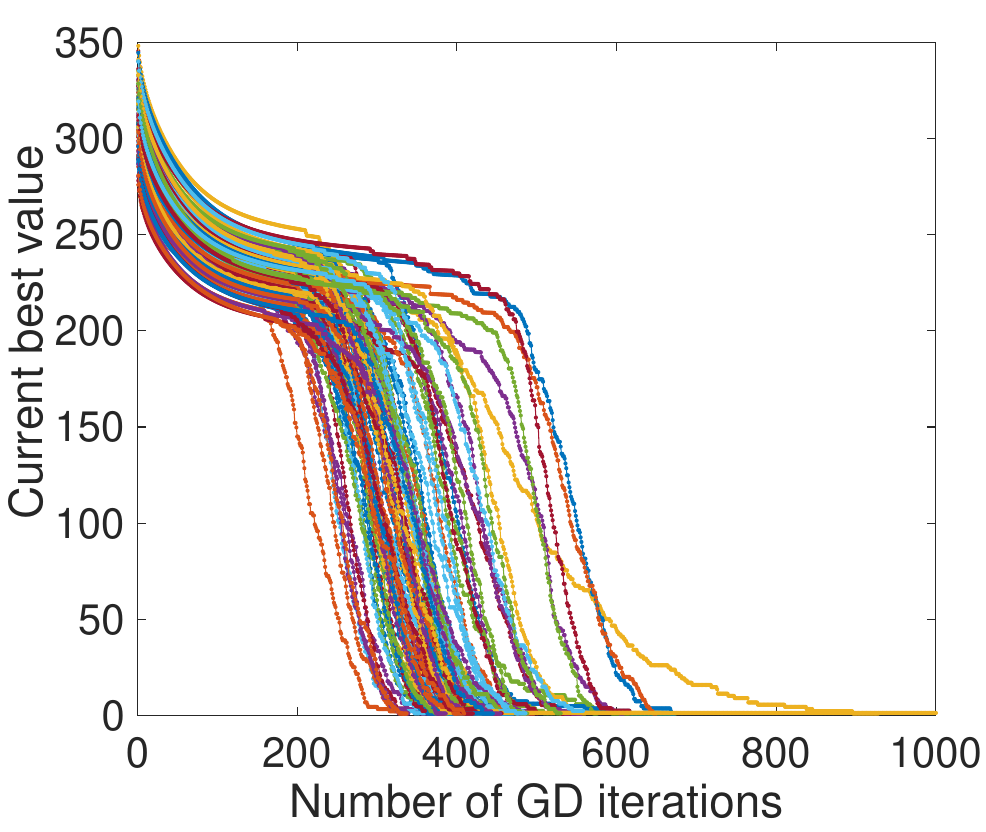}}
	\caption{Left: the number of GD iterations  (in both AM and SVGD) at each time step. 
		Right: the current best value plotted against the GD iterations (in AM) where 
		each line represents a time step.The results are for $\theta =1$ and $M=100$ in the Lorenz 96 example.}
	\label{fig:t1_g1}
\end{figure}

The same sets of figures are also produced for $\theta=0.5$ (Fig.~\ref{fig:t1_b05} for the average bias and Fig.~\ref{fig:t1_g05} for the number of iterations and the current best values) and for $\theta=1$ (Fig.~\ref{fig:t1_b1} for the average bias and Fig.~\ref{fig:t1_g1} for the number of iterations and the current best values). Note that, in Fig.~\ref{fig:t1_b1} the bias of {EnKF} is enormously higher than those of the other methods and so is omitted. The conclusions drawn from these figures are largely the same as those for $\theta=0$, where the key information is that VEnKF significantly outperforms the other methods in terms of estimation bias, and within VEnKF, the results of AM are better than those of SVGD. Regarding the number of GD iterations in AM-VEnKF, one can see that in these two cases (especially in $\theta=1$) it takes evidently more GD iterations for the algorithm to converge,  which we believe is due to the fact that the noise in these two cases are not additive and so the observation models deviate further away from the Gaussian-linear setting. 

\begin{figure}[htbp]
	\centering
	\centerline{
		\includegraphics[width=0.5\textwidth]{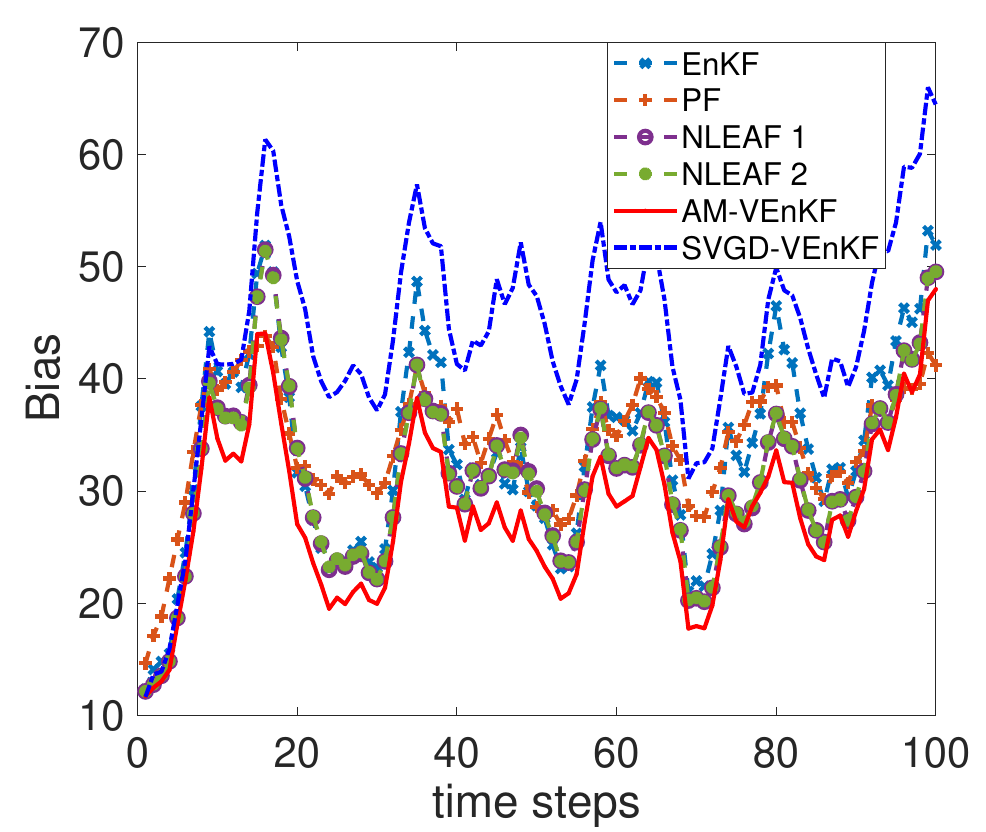}
		\includegraphics[width=0.5\textwidth]{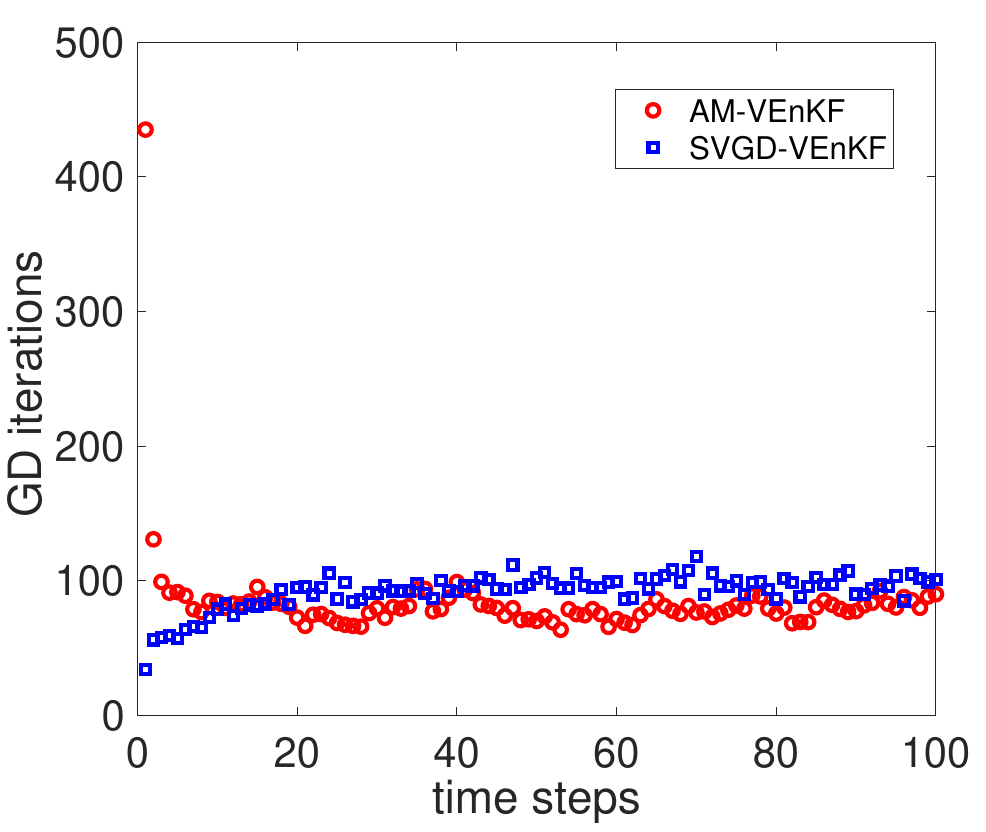}}
	%\kern-\skipamount
	%	\caption{The average bias at each time step for $\epsilon=0$ and $M=20$.}
	%	\label{fig:t2_bl0}
	%\begin{minipage}[htbp]{0.5\linewidth}
	\centerline{\includegraphics[width=0.5\textwidth]{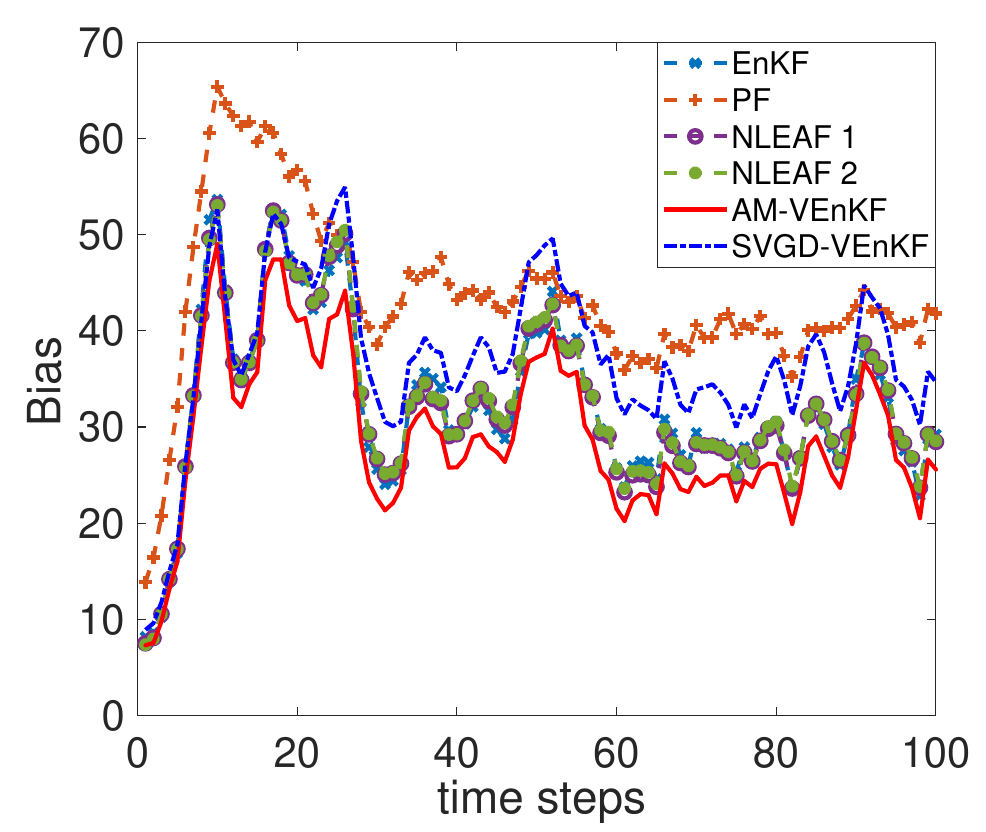}
		\includegraphics[width=0.5\textwidth]{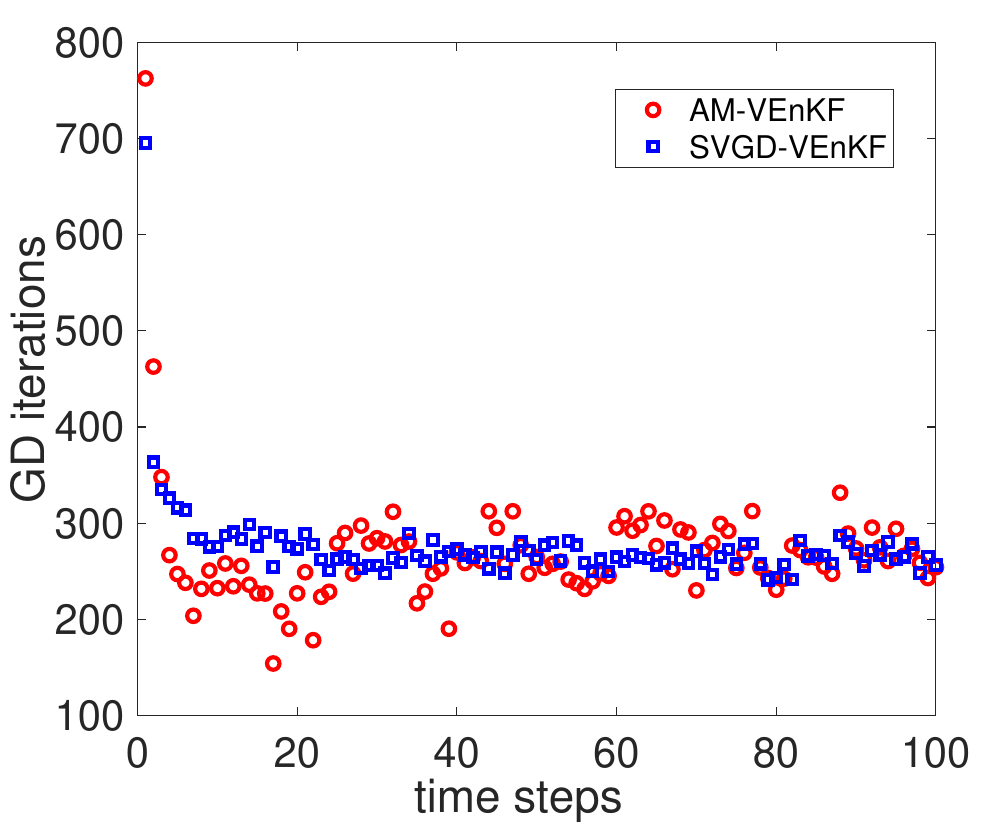}}
	\centerline{\includegraphics[width=0.5\textwidth]{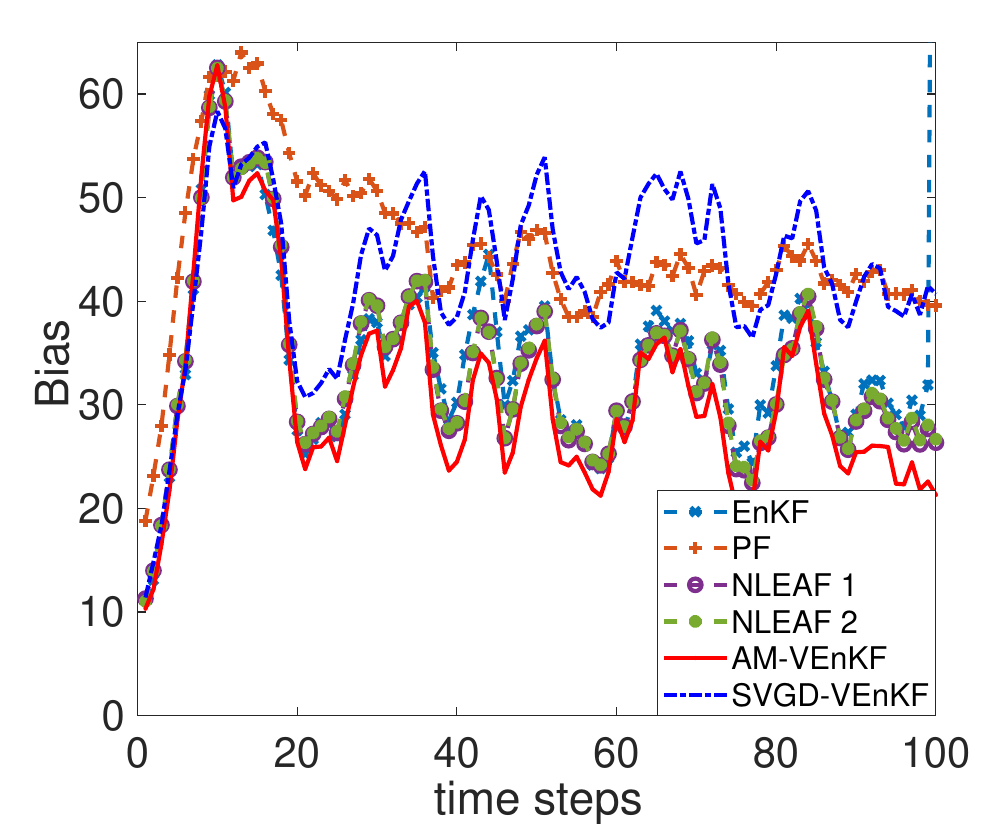}
		\includegraphics[width=0.5\textwidth]{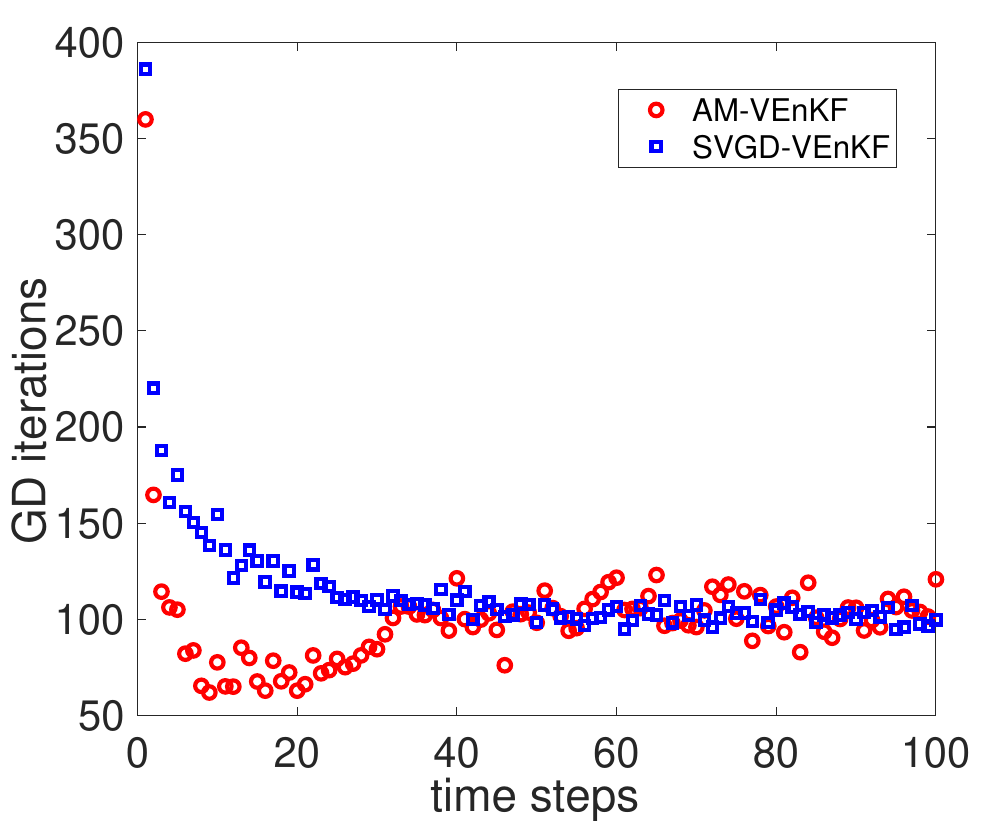}}
	\caption{The results for $M=20$ in the Lorenz 96 example. The figures on the left show the average bias at each time step; the ones on the right show the 
		number of GD iterations (in both AM and SVGD) at each time step.
		From top to bottom are respectively the results of $\theta=0$, 0.5 and 1.}
	\label{fig:t2_bl1}
	%\centering
	%\includegraphics[width=3in]{figs_LMENKF/lor96/iter_lor96_eps1_loc.eps}
	%\caption{t3: iterations(20) , $\theta=0.001$}
\end{figure}

As has been mentioned, we also conduct the experiments for a smaller sample size $M=20$ with localization employed, and we show the average bias results for $\theta=0$, $\theta=0.5$ and $\theta=1$ in Fig.~\ref{fig:t2_bl1}. Similar to the larger sample size case, the bias is also averaged over 200 trials. In this case, we see that the advantage of VEnKF is not as large as that for $M=100$, but nevertheless VEnKF still yields clearly the lowest bias among all the tested methods. On the other hand, the results of the two VEnKF methods are quite similar while that of AM-VEnKF is slightly lower. Also shown in Fig.~\ref{fig:t2_bl1} are the number of GD iterations at each time step for all the three cases, which shows that the numbers of GD iterations used are smaller than their large sample size counterparts. 

\subsection{Fisher's equation}
Our second example is the Fisher's equation, a baseline model of wildfire spreading, where filtering is often needed to assimilate observed data at selected locations into the model~\cite{mandel2008wildland}. Specifically, the Fisher's equation is specified as follows,  
%\numparts
\begin{subequations}
	\begin{eqnarray}
	\label{eq:fisher}
	c_t=Dc_{xx}+rc(1-c),\,\, 0<x<L,\,\, t>0,\\
	c_x(0,t)=0,\,\,c_{x}(L,t)=0,\,\,c(x,0)=f(x),
	\end{eqnarray}
\end{subequations}
%\endnumparts
where $D=0.001$, $r=0.1$, $L=2$ are prescribed constants, and the noise-free initial condition $f(x)$ takes the form of, 
\begin{equation}
f(x)=\left\{\begin{array}{rcl}0, & & 0\leq x<L/4\\
{4x}/{L}-1, & & {L}/{4}\leq x<{L}/{2}\\
3-{4x}/{L}, & & {L}/{2}\leq x<{3L}/{4}\\
0, & & {3L}/{4}\leq x\leq L.
\end{array}\right.
\end{equation}

In the numerical experiments we use an upwind finite difference scheme and discretize the equation onto $N_x= 200$ spatial grid points over the domain $[0,\,L]$, yielding a 200 dimensional filtering problem. The time step size is determined by $D\frac{\Delta t}{\Delta x^2}=0.1$ with $\Delta x=\frac{L}{N_x-1}$ and the total number of time steps is 60. The prior distribution for the initial condition is $U[-5,5]+f(x)$, and in the numerical scheme a model noise is added in each time step and it is assumed to be 
in the form of $N(0,C)$, where 
$$C(i,j)=0.3\exp(-(x_i-x_j)^2/L), \ i,\ j=1,\ldots,N_x,$$ 
with $x_i,x_j$ being the grid points. 

The observation is made at each grid point, and the observation model is as described in Section~\ref{sec:obsmodel}. Once again we test the three cases associated with $\theta=0,\,0.5$ and $1$. The ground truth and the data are both simulated from the model described above. 

\begin{figure}[htbp]
	\centering
	\centerline{\includegraphics[width=0.5\textwidth]{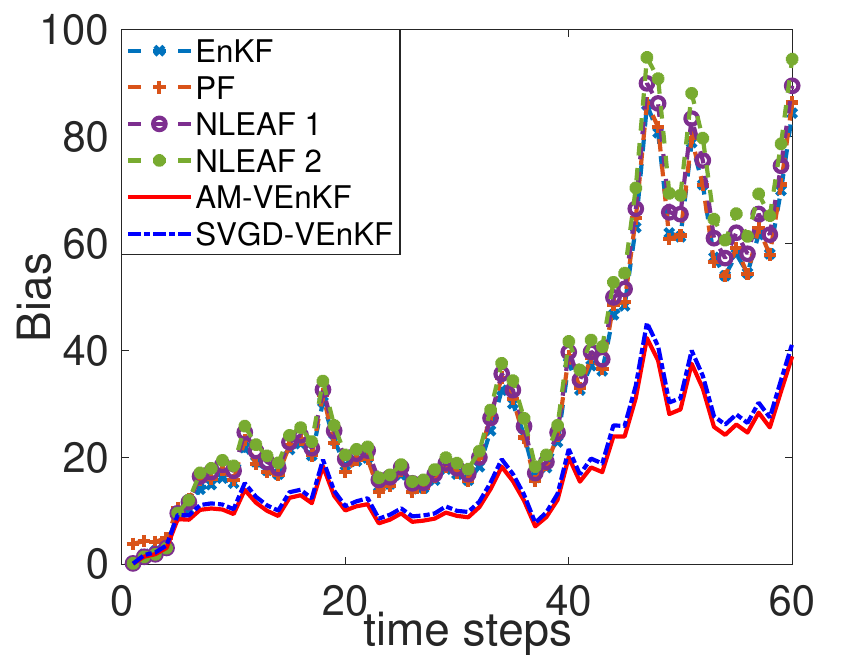}
		\includegraphics[width=0.5\textwidth]{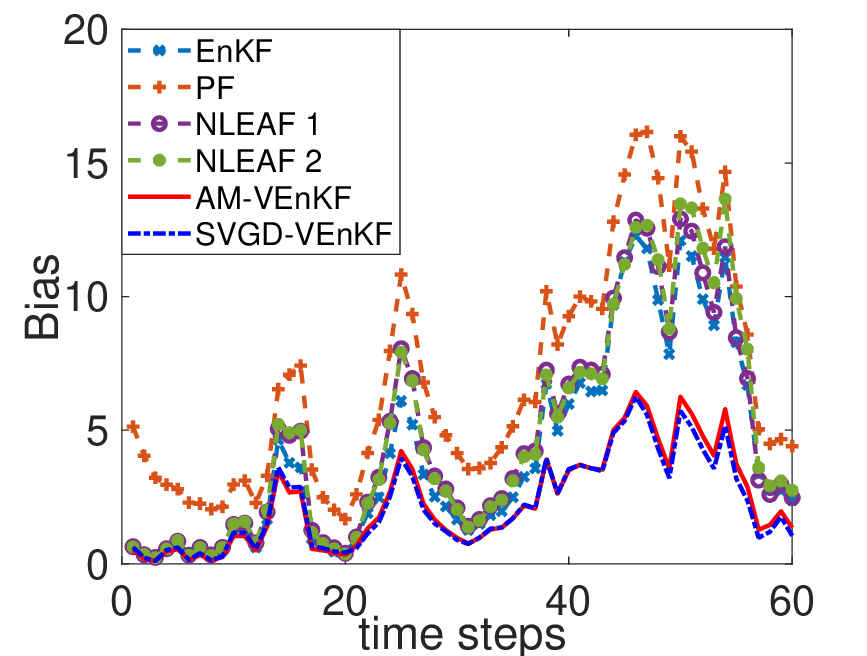}}
	\centerline{\includegraphics[width=0.5\textwidth]{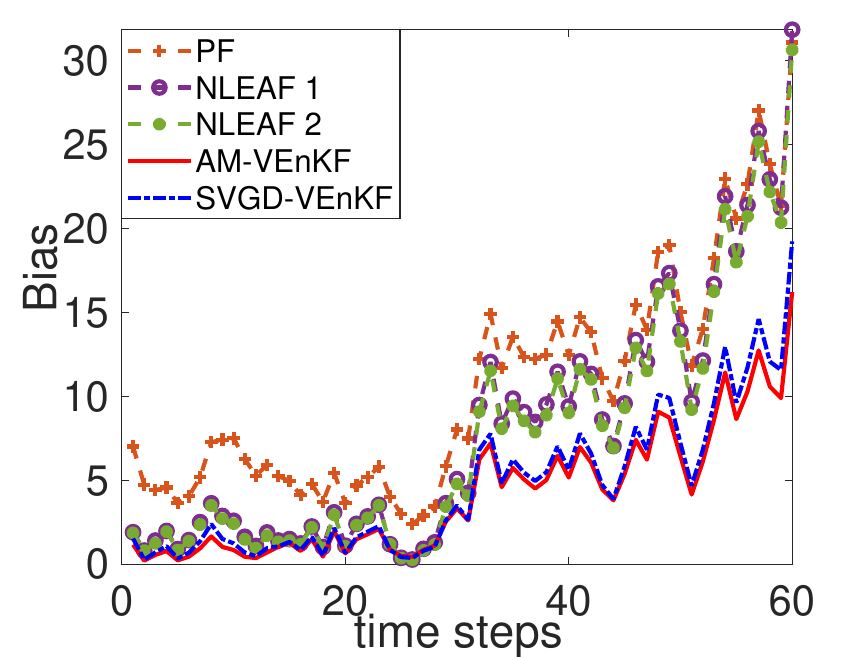}}
	\caption{The average bias at each time step in the Fisher's equation example. From top to bottom: $\theta=0$, $\theta=0.5$ and $\theta=1$.}
	\label{fig:fisher}
	%\centering
	%\includegraphics[width=3in]{figs_LMENKF/lor96/iter_lor96_eps1_loc.eps}
	%\caption{t3: iterations(20) , $\theta=0.001$}
\end{figure}

We test the same set of filtering methods as those in the first example. Since in practice, it is usually of more interest to consider a small ensemble size relative to the dimensionality, we choose to use 50 particles for this 200 dimensional example. Since the sample size is smaller than the dimensionality, the sliding window localization with $l=5$ and $k=3$ is used. All the simulations are repeated 200 times and the average biases are plotted in Fig.~\ref{fig:fisher} for all the three cases ($\theta=0,\,0.5$ and $1$). We see that in all the three cases the two VEnKF methods result  in the lowest estimation bias among all the methods tested, and the results of the two VEnKF methods are rather similar. It should be mentioned that, in the case of $\theta=1$, the bias of {EnKF} is omitted as it is  enormously higher than those of the other methods.

\begin{figure}[htbp]
	\centerline{\includegraphics[width=0.33\textwidth]{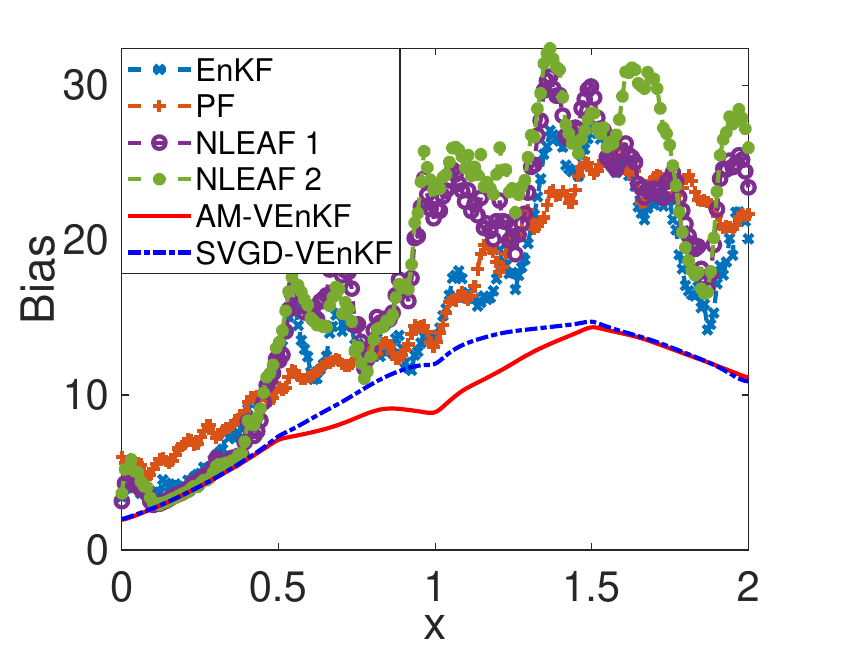}
		\includegraphics[width=0.33\textwidth]{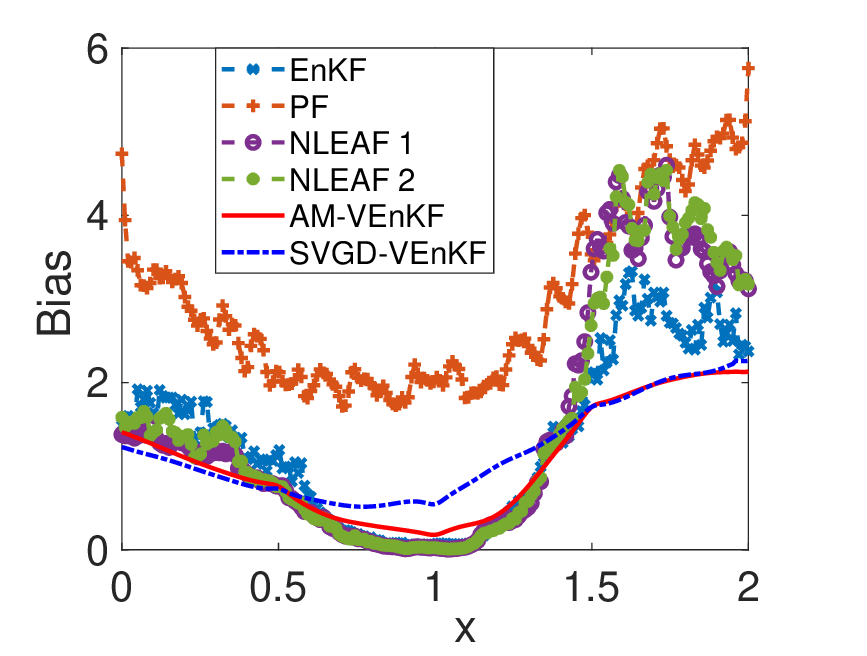}
		\includegraphics[width=0.33\textwidth]{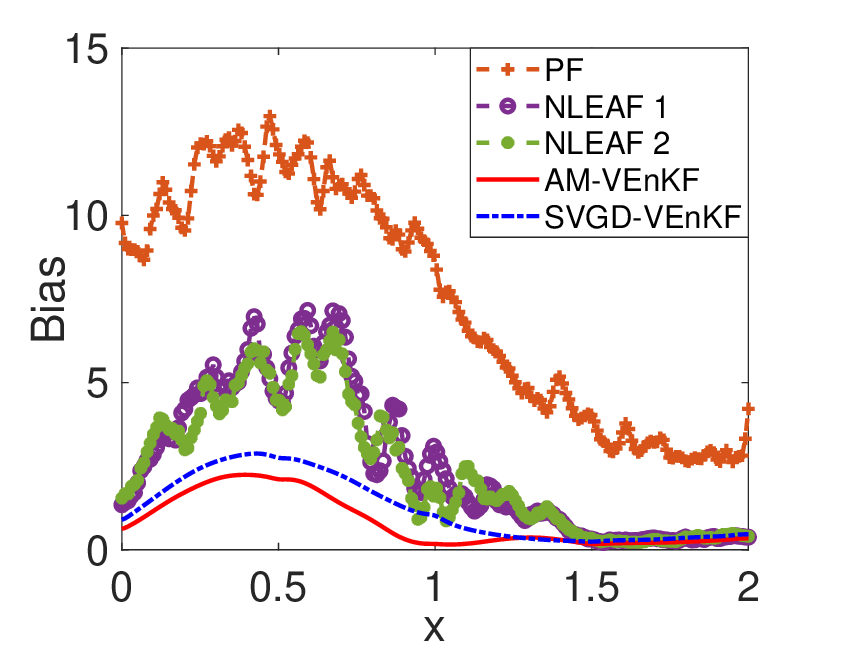}}
	\centerline{\includegraphics[width=0.33\textwidth]{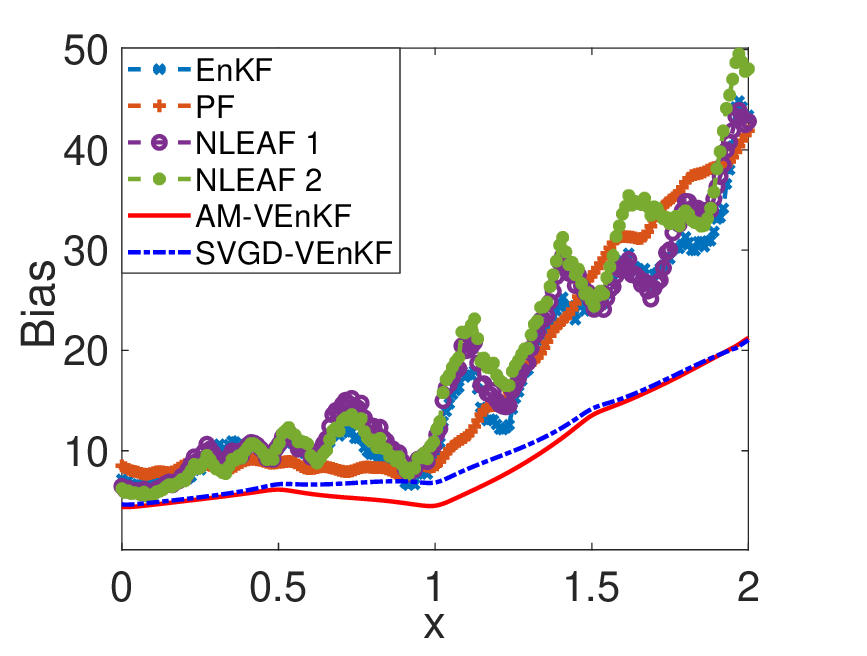}
		\includegraphics[width=0.33\textwidth]{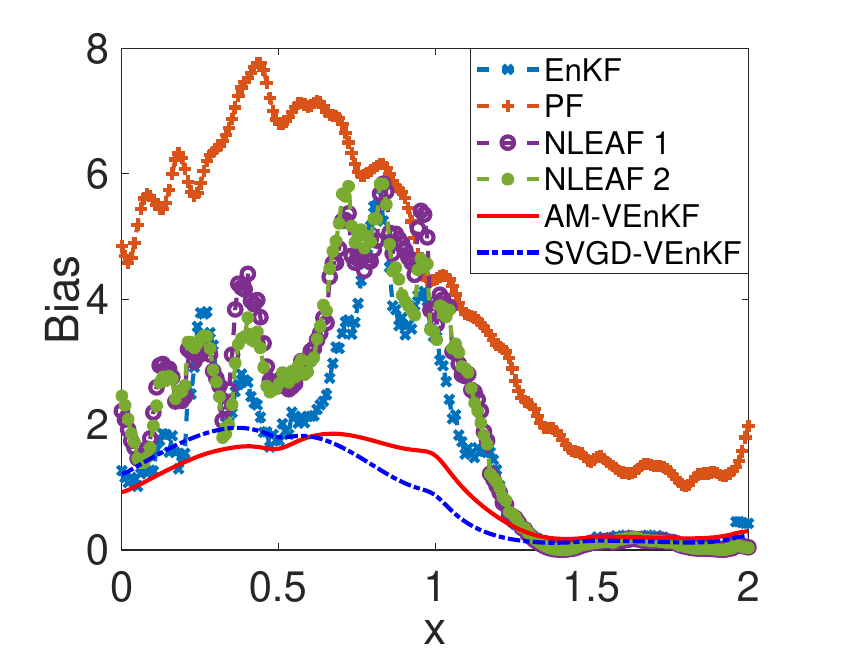}
		\includegraphics[width=0.33\textwidth]{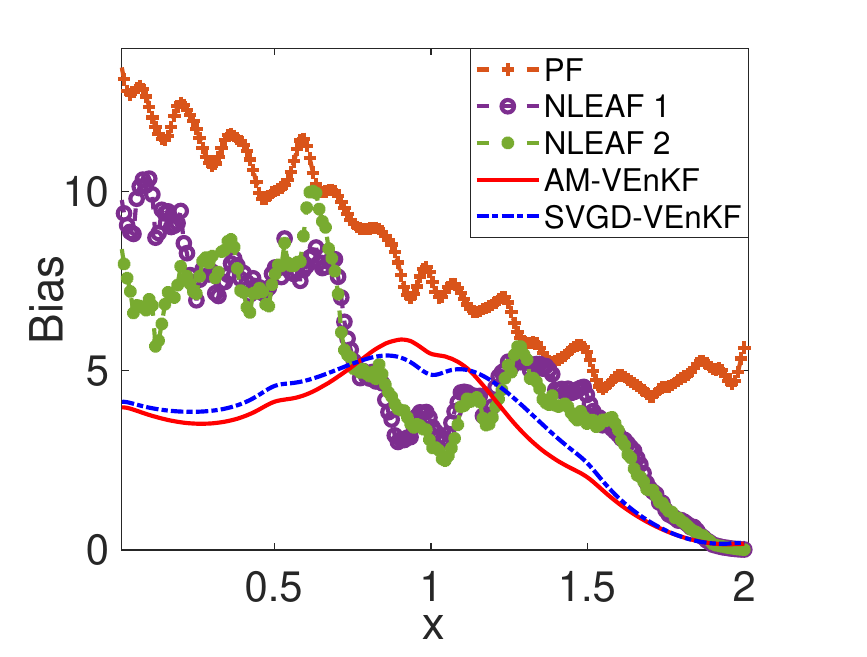}}
	\centerline{\includegraphics[width=0.33\textwidth]{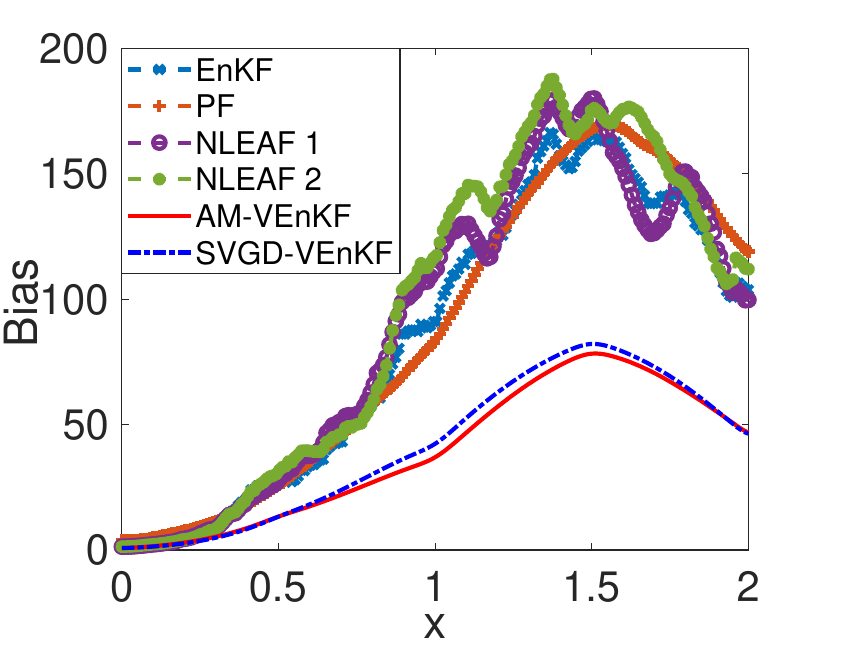}
		\includegraphics[width=0.33\textwidth]{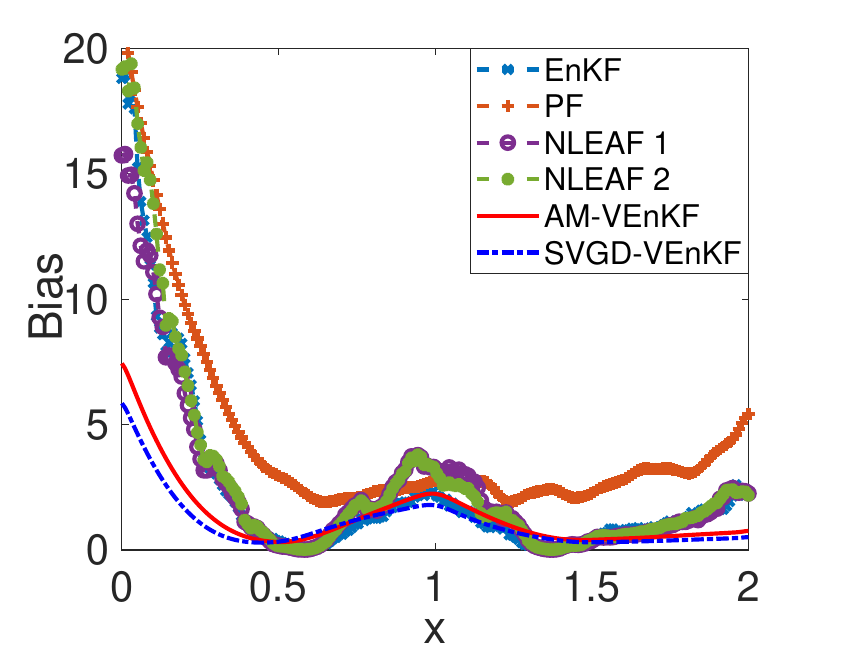}
		\includegraphics[width=0.33\textwidth]{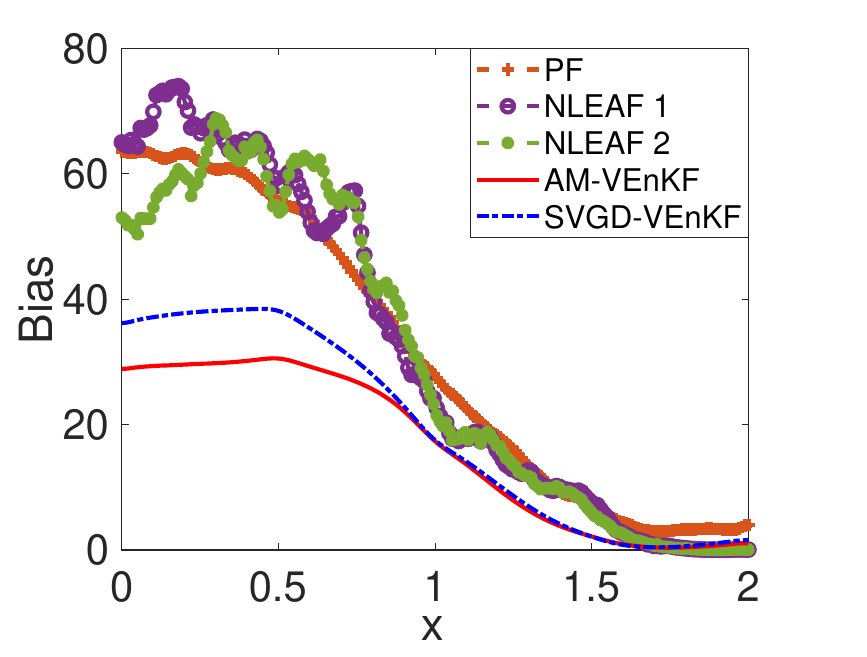}}
	\caption{The estimation bias at $t=10$ (top), $t=30$ (middle) and $t=60$ (bottom), in the Fisher's equation example. From left to right: $\theta=0$, $\theta=0.5$ and $\theta=1$.}
	\label{fig:fisher2}
	%\centering
	%\includegraphics[width=3in]{figs_LMENKF/lor96/iter_lor96_eps1_loc.eps}
	%\caption{t3: iterations(20) , $\theta=0.001$}
\end{figure}

As the bias results shown in Fig.~\ref{fig:fisher} are averaged over all the dimensions, it is also useful to examine the bias at each dimension. We therefore plot in Fig.~\ref{fig:fisher2} the bias of each grid point at three selected time steps $t=10,\,30,$ and 60. The figures illustrate that, at all these time steps, the VEnKF methods yield substantially lower bias at the majority of the grid points, which is consistent with the average bias results shown in Fig.~\ref{fig:fisher}. 
%We have also examined other time steps where the results are similar and therefore not presented here. 
We also report that, the wall-clock time for solving the optimization problem in each time step in AM-VEnKF is approximately 2.0 seconds (on a personal computer with a 3.6GHz processor and 16GB RAM), indicating a modest computational cost  in this 200 dimensional example. 

\subsection{Lorenz 2005 model}
Here we consider the Lorenz 2005 model~\cite{lorenz2005designing} which products spatially more smoothed model trajectory than Lorenz 96. The Lorenz 2005 model is written in the following scheme,
\begin{equation}\label{eq:lorII}
\frac{dx^n}{dt}=[x,x]^{K,n}-X^n+F, \quad n=1,\ldots,N.
\end{equation}
where \[[x,x]^{K,n}=\sum\limits_{j=-J}^{J}{'}\sum\limits_{i=-J}^{J}{'}(-x^{n-2K-i}x^{n-K-j}+x^{n-K+j-i}x^{n+K+j})/K^2,\] and this equation is composed with periodic boundary condition. $F$ is the forcing term and $K$ is the smoothing parameter while $K<<N$, and one usually sets $J=\frac{K-1}{2}$ if $K$ is odd, and $J=\frac{K}{2}$ if $K$ is even. Noted that the symbol $\sum{'}$ denote a modified summation which is similarly with generally summation $\sum$ but the first and last term are divided by $2$. Moreover if $K$ is even the summation is $\sum{'}$, and if $K$ is odd the summation is replaced by ordinary $\sum$. 

It is worth noting that, when setting $K=1$, $N=40$, and $F=8$, the model reduces to Lorenz 96. In this example, we set the model as $N=560$, $F=10$ and $K=16$,  resulting in a 560-dimensional filtering problem. Following the notations in Sec.~\ref{sec:lor96}, Lorenz 2005 is also represented by a standard discrete-time fourth-order Runge-Kutta solution of Eq.~(\ref{eq:lorII}) with $\Delta t=0.01$ where the same model noise is added, and the state and observation pair $\{{\bf{x}}_t,{\bf{y}}_t\}$ is similarly denoted by Eq.~(\ref{eq:discretelorenz}). We reinstate that in this example the observation model is chosen differently (see Sec.~\ref{sec:obsmodel}).
% \begin{eqnarray}
% \bf{x}_t &=& F(\bf{x}_{t-1})+\alpha_t,\nonumber\\
% \bf{y}_t &=& M(\alphax_t)+M(\alphax_t)^{\theta}\beta_t,\quad t=1,2,\ldots
% \end{eqnarray}
%While in this example we set another observed operator $M(\bf{x}_t)=\exp(\alpha \bf{x}_t)$ and $\alpha$ is the adjusted nonlinear parameter in observation operator $M$, here is $0.5$. 
And the initial state is chosen to be ${\bf{x}}_0\sim U[0,5]$.
\begin{figure}[htbp]
	\centerline{
		\includegraphics[width=0.5\textwidth]{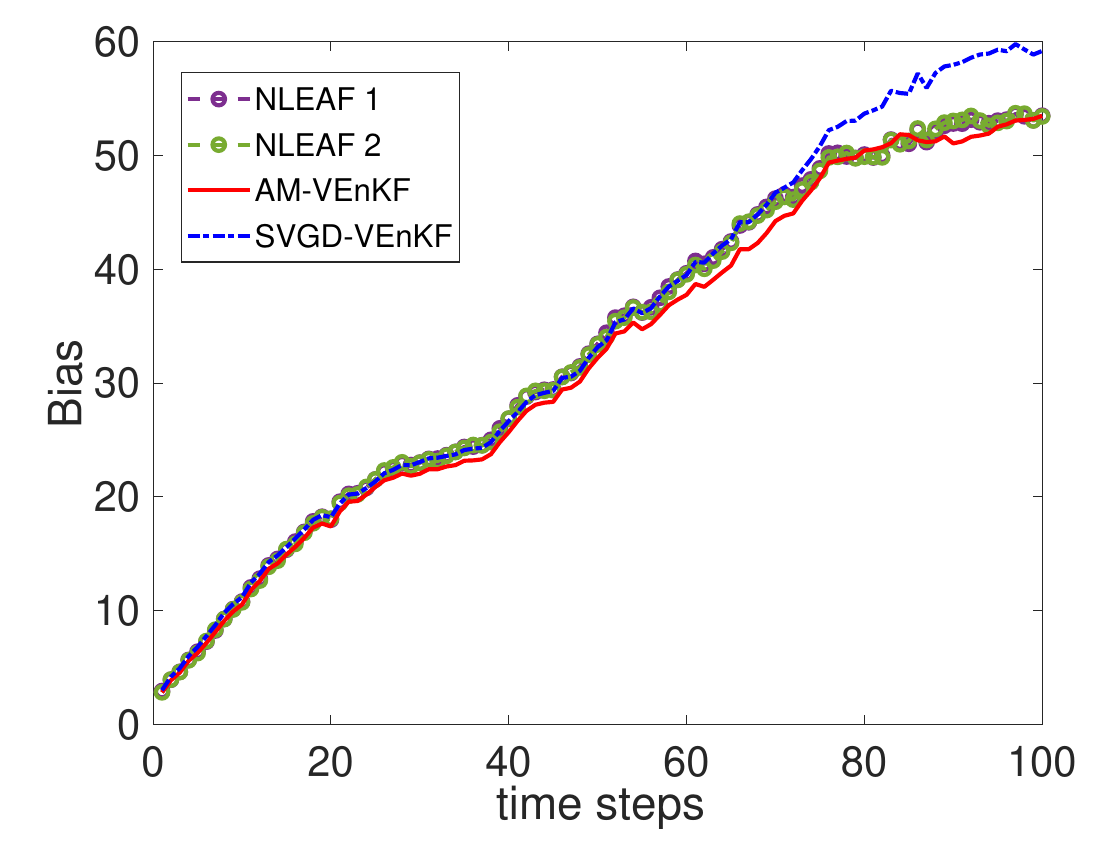}
		\includegraphics[width=0.5\textwidth]{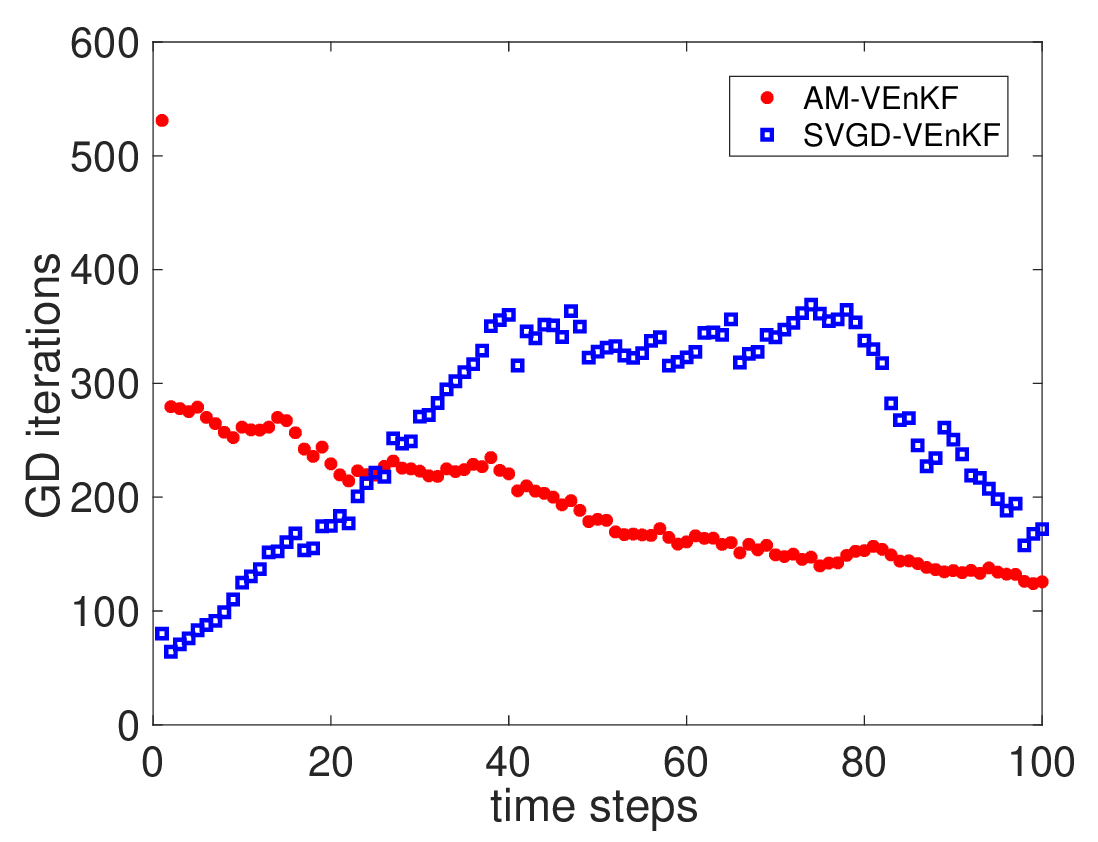}}
	\centerline{
		\includegraphics[width=0.5\textwidth]{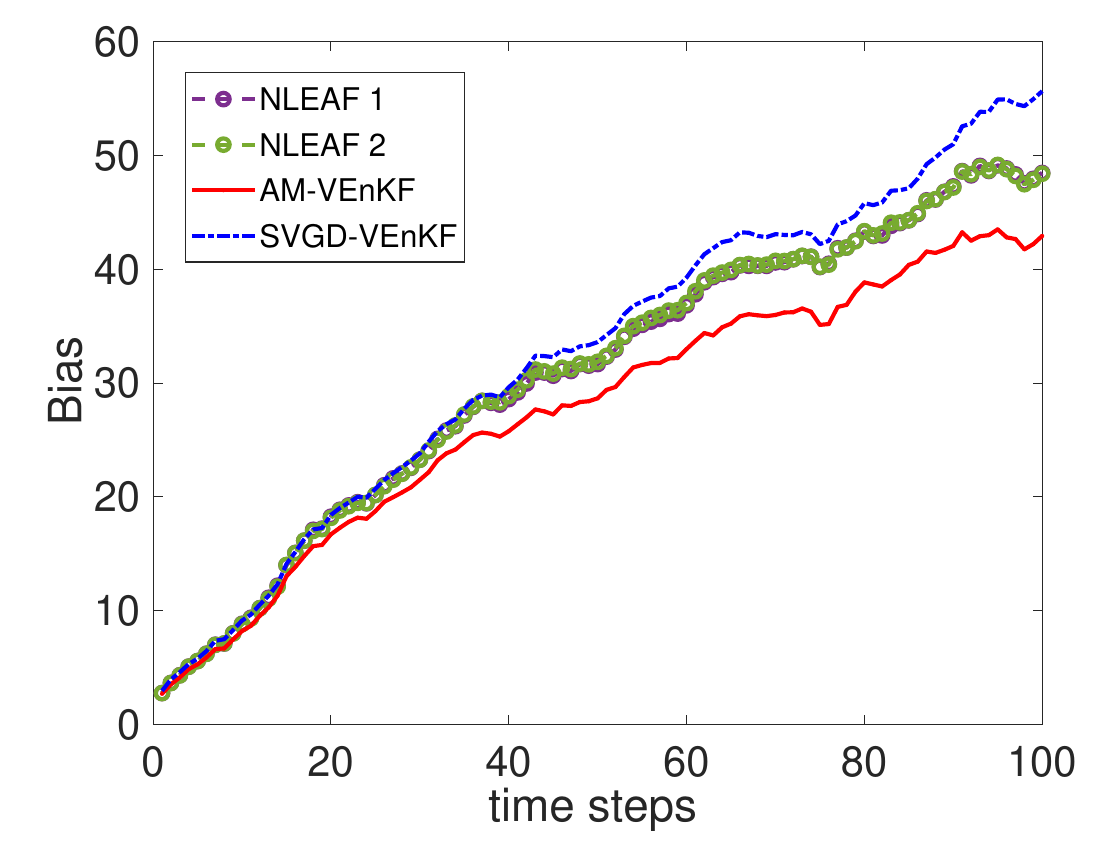}
		\includegraphics[width=0.5\textwidth]{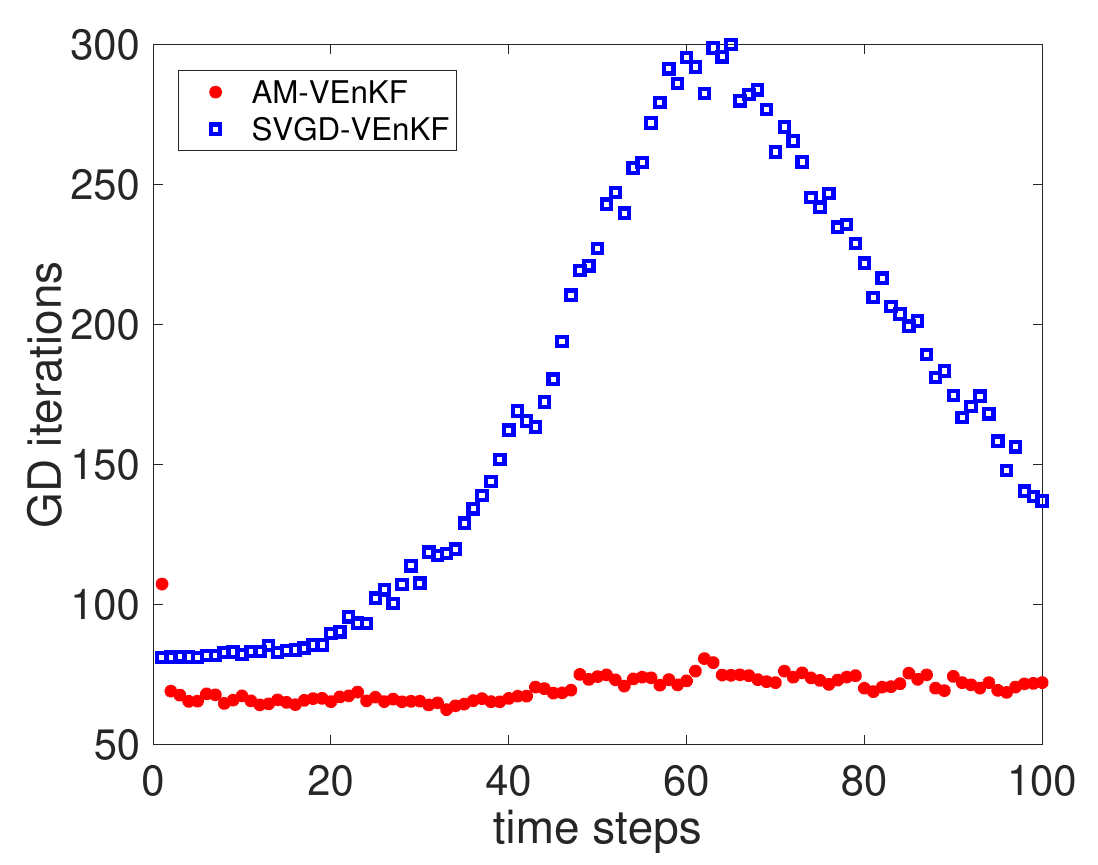}}
	\centerline{
		\includegraphics[width=0.5\textwidth]{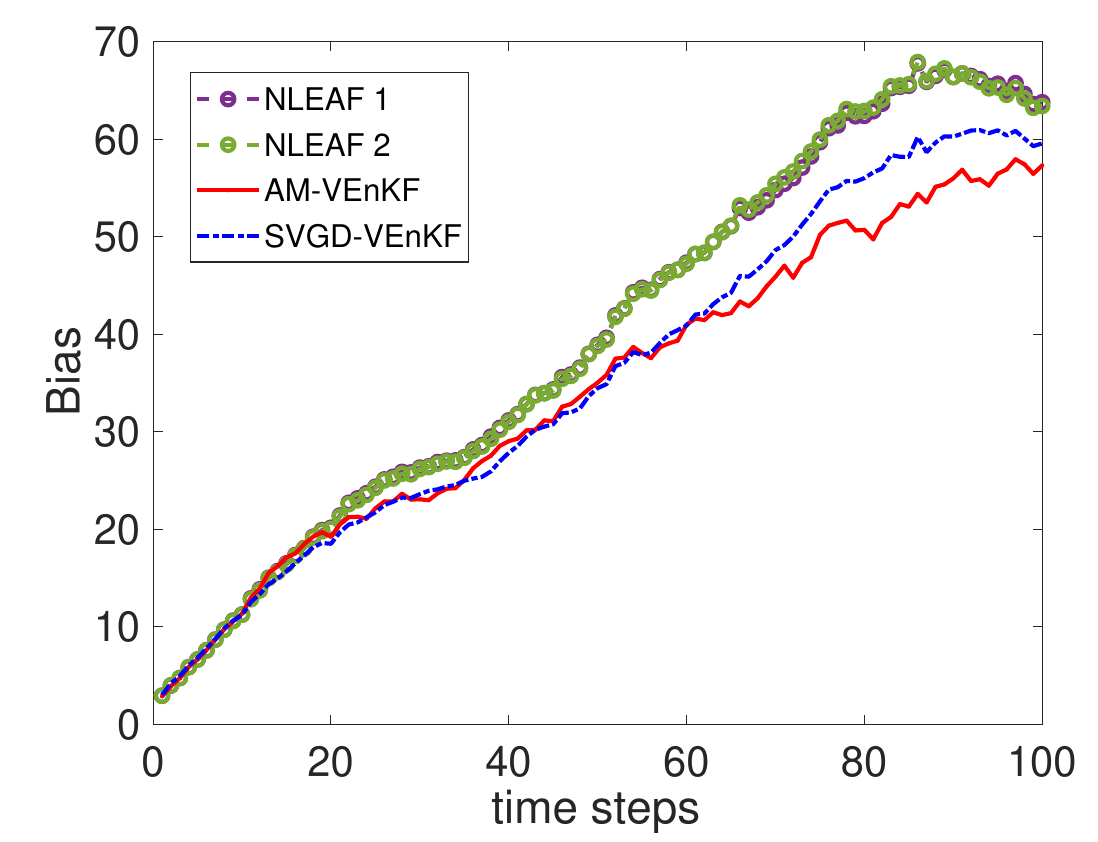}
		\includegraphics[width=0.5\textwidth]{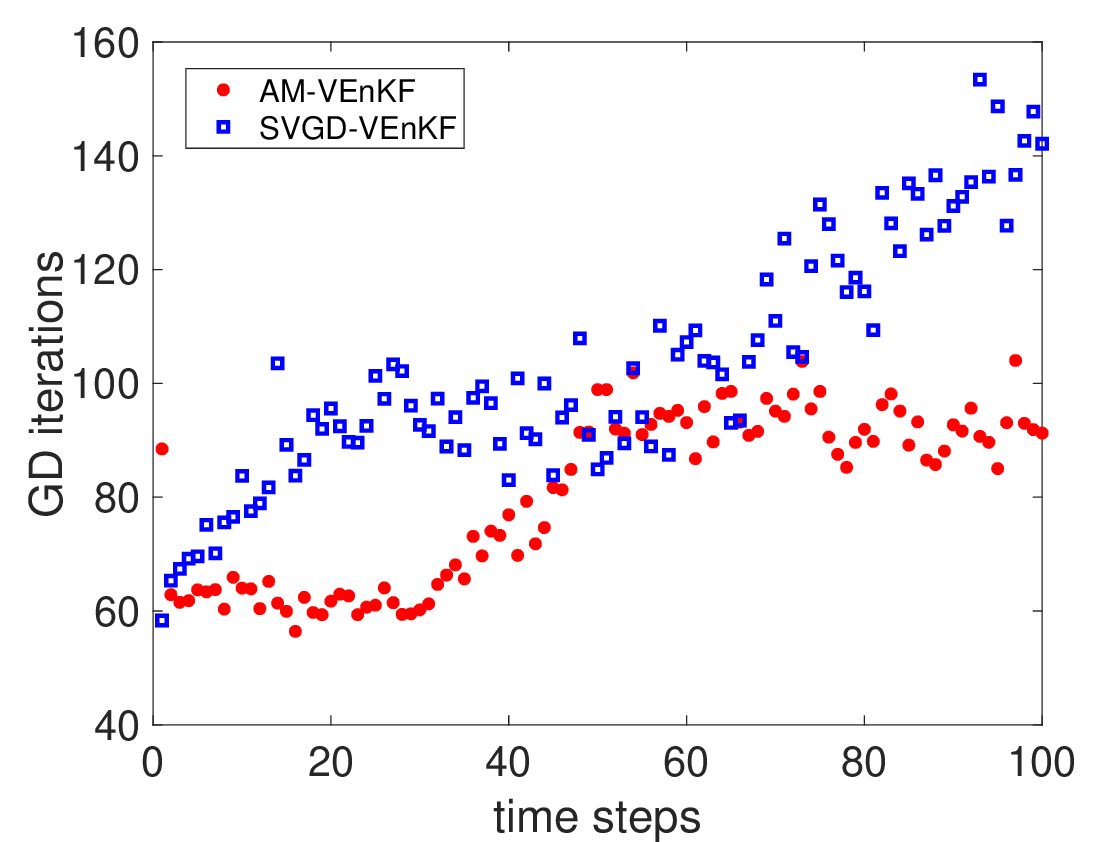}}
	%\kern-\bigskipamount
	\caption{The results for the Lorenz 2005 example: the figures on the left show the average bias at each time step; the ones on the right show the number of GD iterations (in both AM and SVGD) at each time step.
		From top to bottom are respectively the results of $\theta=0$, 0.5 and 1.}
	\label{fig:t3}
\end{figure}

In this numerical experiments, we test the same set of methods as those in the first two examples, where in each method 100 particles are used. Due to the small ensemble size, it is necessary to adopt the sliding-window localization with 
$(l,k)=(5,3)$ in all methods except PF. We observe that the errors in the results of EnKF and PF are significantly larger than those in the other methods, and so those results are not presented here. It should be noted that the stopping threshold is as $\Delta_F=0.5$ during nearest $\Delta_k=20$ iterations in AM-VEnKF. All methods are repeated 20 times and we plot the averaged bias and the averaged GD iterations for all the three cases ($\epsilon=0$, $0.5$ and $1$) in Fig.~\ref{fig:t3}. One can see from the figures that, in the first case ($\epsilon=0$) the results of all the methods are quite similar, while in the other two cases, the results of AM-VEnKF are clearly better than those of all the other methods. 

\section{Closing Remarks}\label{sec:closing}

%\kern-\bigskipamount
%In summary we present a affine mapping based variational ensemble Kalman filter, which can be applied to generic 
%observation models. The method uses a variational formulation and the resulting KLD minimization problem is solved with a gradient descent scheme. 
We conclude the paper with the following remarks on the proposed VEnKF framework. First 
%we reinstate that, the KLD minimization step in LMEKF does not involve the underlying dynamical model which is usually computationally intensive.
we reinstate that, the Fisher's equation example demonstrates that the KLD minimization problem in AM-VEnKF can be solved rather efficiently, and more importantly this optimization step does not involve simulating the underlying dynamical model. As a result, this step, though more complicated than the update in the standard EnKF, may not be the main contributor to the total computational burden, especially when the underlying dynamical model is computational intensive. Second, it is important to note that, although VEnKF can deal with generic observation models, it still requires that the posterior distributions are reasonably close to Gaussian, an assumption needed for all EnKF type of methods. For strongly non-Gaussian posteriors, it is of our interest to explore the possibility of incorporating VEnKF with some existing extensions of {EnKF} that can handle strong non-Gaussianity, such as the mixture Kalman filter~\cite{stordal2011bridging}. Finally, in this work we provide two transform mappings, the affine mapping and the RKHS mapping in the SVGD framework. {In the numerical examples studied here, the affine mapping seems to achieve a better performance, but we acknowledge that more comprehensive comparisons should be done to understand the advantages and limitations of different types of mappings. A related issue is that, some existing works such as~\cite{pulido2019sequential} use more flexible and complicated mappings and so that they can approximate arbitrary posterior distributions. It is worth noting, however, this type of methods are generally designed for problems where a rather large number of particles can be afforded, and therefore are not suitable for the problems considered here. Nevertheless, developing more flexible mapping based filters is an important topic that we plan to investigate in future studies. }

%\section*{Acknowledgment}
%
%
%The work was partially supported by NSFC under grant ***.

% or


\appendix

\section{SVGD-VEnKF}

In this section, we discuss the procedure for constructing the mapping using the Stein variational gradient descent (SVGD) formulation \cite{liu2016stein}, which provides a nonlinear transform from the prior to the posterior in each time step.
%In what follows we first describe SVGD in a generic setting and then discuss how it is used to construct the mapping in our framework. 

%%SVGD is a nonparametric variational inference algorithm, which iteratively transform a set of particles from prior distribution to  approximate a given %posterior distribution via a series of mappings. Here we give a quick overview of its main idea in next subsection.
%
%\subsection{An introduction to SVGD}
%
%Suppose we want to determine a stochastic process $X$ in $\mathbb{R}^{d_x}$, however we only know the data which are from another process $Y$ in $\mathbb{R}^{d_y}$. The relationship between the processes is given by a known nonlinear observation system such that
%\begin{equation}
%\label{eq:obe}
%Y=G(X,\eta)
%\end{equation}
%where $\eta$ is the random observational error, and we can get the likelihood distribution $\pi(Y|X)$ by Eq.~(\ref{eq:obe}). This is known as Bayesian Inverse problem, which use observation data $Y$ to infer the real state $X$. And  we assume that the prior density $\pi(X)$ is known and can be expressed by some distributions.

Recall that in Section~\ref{sec:lmekf} we want to find a mapping by solving   
\begin{equation}
\min_{T\in \mathcal{H}}\kld(\pi_T,q),\label{e:kld_svgd}
\end{equation}
where $q(\cdot)=\hat{\pi}(\cdot|y_{1:t})$ and $\mathcal{H}$ is a certain function space that will be specified later. 

Following the same argument in Sec.~\ref{sec:lof}, we obtain that Eq.~(\ref{e:kld_svgd}) is equivalent to, 
\begin{equation}
\min_{T\in \mathcal{H}}\kld(p(\tilde{x}),q_{T^{-1}}(\tilde{x})),\label{e:svgdminkld2}
\end{equation}
where $q_{T^{-1}}(\cdot)$ is as defined in Section~\ref{sec:lof}. 

Now we need to determine the function space $\mathcal{H}$. While in the proposed AM-VEnKF method $\mathcal{H}$ is chosen to be an affine mapping space, the SVGD framework specifies $\mathcal{H}$ via a reproducing kernel Hilbert space (RKHS)~\cite{scholkopf2018learning}.

First we write the mapping  $T$ in the form of,
\begin{equation}
\label{eq:mp_svgd}
T(\tilde{x})=\tilde{x}+\tau \phi(\tilde{x}),
\end{equation}
where $\tau$ is a prescribed stepsize. 
% the variable $x$ express the transformed variable, and $\epsilon$ is the stepsize in mapping, while $\pi_T(x)=|\nabla T^{-1}|\pi(\tilde{x})$. 
Next we assume that mapping $\phi$ is chosen from a RKHS $\mathcal{H}_K$ specified by a reproducing kernel $K(\cdot,\cdot)$. Therefore the optimisation problem~(\ref{e:svgdminkld2}) becomes,
\begin{equation}
\label{eq:mkl_svgd}
\min_{\phi\in\mathcal{H}_K}\kld(p(\tilde{x}),q_{T^{-1}}(\tilde{x})).
\end{equation}
%
% The Gateaux derivative of the Kullback-Leibler divergence, in the direction $\phi$ is given by
%\begin{equation}
%	\label{eq:dkld_svgd}
%	\nabla_{\phi}\kld(p(\tilde{x}),q_{T^{-1}}(\tilde{x}))=-\int d_{\tau}\log[q_{T^{-1}}(\tilde{x})]|_{\tau=0}p(\tilde{x})d\tilde{x},
%\end{equation}
%where the derivative of the transformed log-posterior density is
%\[d_{\tau}\log[q_{T^{-1}}(\tilde{x})]=\nabla_{T}\log[q(T(\tilde{x}))]'d_{\tau}T+\Tr(\nabla_{\tilde{x}}T^{-1}d_{\tau}\nabla_{\tilde{x}}T).\]
%
%Rewrite the Eq.~(\ref{eq:dkld_svgd}) using the mapping Eq.~(\ref{eq:mp_svgd}), then the directional derivative of KL results in,
%\begin{equation}
%	\nabla_{\phi}{\kld(p(\tilde{x}),q_{T^{-1}}(\tilde{x}))}=-\int[\nabla_{\tilde{x}}\log[q(\tilde{x})]'\phi(\tilde{x})+\Tr(\nabla\phi)]p(\tilde{x})d{\tilde{x}}.
%\end{equation}
%
%Taking the derivative of KL divergence w.r.t the perturbation magnitude $\tau$, so theoptimization problem is taken to be
%\begin{eqnarray}
%	\phi^*=\mathrm{E}_{\tilde{x}\sim p}[\nabla_{\tilde{x}}\log[q(\tilde{x})]'\phi(\tilde{x})+\Tr(\nabla_{\tilde{x}}\phi(\tilde{x})).]
%\end{eqnarray}
%Where in the term of mapping operator $\mathcal{A}'\phi(\tilde{x})=\nabla_{\tilde{x}}\log[q(\tilde{x})]'\phi(\tilde{x})+\Tr(\nabla_{\tilde{x}}\phi(\tilde{x}))$, $\mathcal{A}$ is called the Stein operator.
% 
In the SVGD framework, one does not seek to solve the optimisation problem in Eq.~(\ref{eq:mkl_svgd}) directly; instead it can be derived that the direction of steepest
descent is
%the mapping function space $\mathcal{H}$ is defined in RKHS, so using $\phi(x)=\langle K(\cdot,x),\phi(x)\rangle _{\mathcal{H}}$, than the optimal mapping function $\phi^*$ in filtering time $t$ can have a simple closed-form solution, that is
\begin{equation}
\phi^*(\cdot)=\mathrm{E}_{\tilde{x}\sim p}[\nabla_{\tilde{x}}\log q(\tilde{x})K(\tilde{x},\cdot)+\nabla_{\tilde{x}}K(\tilde{x},\cdot)]. \label{e:phiopt}
\end{equation}
It should be noted that we omit the detailed derivation of Eq.~(\ref{e:phiopt})  here and interested readers may consult \cite{liu2016stein} for such details. 
The obtained mapping $\phi^*$ is then applied to the samples which pushes them toward the target distribution. This procedure is repeated until certain stopping conditions are satisfied. The complete SVGD based VEnKF algorithm is given in Alg.~\ref{alg:SVGD}. 
Finally we note that, in the numerical experiments we use the squared exponential kernel with bandwidth $h$:
\[K(x,x') =\exp (-\|x-x'\|_2^2/h),\]
where the implementation details can be found in \cite{liu2016stein}.
%The algorithm using SVGD update based VEnKF in recursive Bayesian filtering is summarized in Alg.~\ref{alg:SVGD}.
\begin{algorithm}
	\caption{SVGD based variational EnKF (SVGD-VEnKF)} %
	\label{alg:SVGD}
	\begin{itemize}
		\item Prediction: 
		\begin{itemize}
			\item Let $\tilde{x}_t^m= F_t(x_{t-1}^m,\alpha_t^m), m=1,2,\ldots,M$;
			\item Let $\hat{\pi}(\cdot|y_{1:t-1}) = N(\tilde{\mu}_t,\tilde{\Sigma}_t)$ where $\tilde{\mu}_t$ and $\hat{\Sigma}_t$
			are computed using Eq.~(\ref{eq:emp});
		\end{itemize}
		\item Update:
		\begin{itemize}
			\item Let $q(\cdot)=\hat{\pi}(\cdot|y_{1:t})\propto \hat{\pi}(\cdot|y_{1:t-1})\pi(y_t|\cdot)$;
			\item Repeat the following steps until the stopping conditions are satisfied;
			\begin{itemize}
				\item Let
				\begin{equation}
				\tilde{\phi}^*(\cdot)= \frac{1}{M}\sum\limits_{m=1}^M[\nabla_{\tilde{x}_t^m}\log q(\tilde{x}_t^m)K(\tilde{x}_t^m,\cdot)+\nabla_{\tilde{x}_t^m}K(\tilde{x}_t^m,\cdot)].\nonumber
				\end{equation}
				\item Let $\tilde{x}_t^m\leftarrow\tilde{x}_t^m+\tau\tilde{\phi}^*({\tilde{x}_t^m}$), $m=1,\ldots,M$. 
			\end{itemize}		
			\item Let
			$x_t^m=\tilde{x}_t^m$, for $m=1,\ldots,M$.
		\end{itemize}
	\end{itemize}
\end{algorithm}

%% you can choose not to have a title for an appendix
%% if you want by leaving the argument blank
%\section{}


%\ifCLASSOPTIONcaptionsoff
%  \newpage
%\fi

%\bibliographystyle{cas-model2-names}
\bibliographystyle{plain}
\bibliography{svgd.bib}

\end{document}